\documentclass{amsart}
\usepackage{amsmath}
\usepackage{amssymb}
\usepackage{amsfonts}

\newtheorem{theorem}{Theorem}
\theoremstyle{plain}

\newtheorem{corollary}{Corollary}

\newtheorem{definition}{Definition}

\newtheorem{lemma}{Lemma}

\newtheorem{remark}{Remark}

\numberwithin{equation}{section}
\input{tcilatex}

\begin{document}
\title[Parabolic sublinear operators with rough kernel]{Parabolic sublinear
operators with rough kernel generated by parabolic Calder\'{o}n-Zygmund
operators and parabolic local Campanato space estimates for their
commutators on the parabolic generalized local Morrey spaces }
\author{F.GURBUZ}
\address{ANKARA UNIVERSITY, FACULTY OF SCIENCE, DEPARTMENT OF MATHEMATICS,
TANDO\u{G}AN 06100, ANKARA, TURKEY }
\email{feritgurbuz84@hotmail.com}
\urladdr{}
\thanks{}
\curraddr{ }
\urladdr{}
\thanks{}
\date{}
\dedicatory{}
\thanks{}

\begin{abstract}
In this paper, the author introduces parabolic generalized local Morrey
spaces and gets the boundedness of a large class of parabolic rough
operators on them. The author also establihes the {parabolic local Campanato}
space estimates for their commutators on parabolic generalized local Morrey
spaces. As its special cases, the corresponding results of parabolic
sublinear operators with rough kernel and their commutators can be deduced,
respectively. At last, parabolic Marcinkiewicz operator which satisfies the
conditions of these theorems can be considered as an example.
\end{abstract}

\maketitle

\section{Introduction}

Let ${\mathbb{R}^{n}}$ be the $n-$dimensional Euclidean space of points $%
x=(x_{1},...,x_{n})$ with norm $|x|=\left( \dsum
\limits_{i=1}^{n}x_{i}^{2}\right) ^{\frac{1}{2}}$. Let $B=B(x_{0},r_{B}) $
denote the ball with the center $x_{0}$ and radius $r_{B}$. For a given
measurable set $E$, we also denote the Lebesgue measure of $E$ by $%
\left
\vert E\right \vert $. For any given $\Omega \subseteq {\mathbb{R}^{n}%
}$ and $0<p<\infty $, denote by $L_{p}\left( \Omega \right) $ the spaces of
all functions $f$ satisfying%
\begin{equation*}
\left \Vert f\right \Vert _{L_{p}\left( \Omega \right) }=\left( \dint
\limits_{\Omega }\left \vert f\left( x\right) \right \vert ^{p}dx\right) ^{%
\frac{1}{p}}<\infty .
\end{equation*}%
Let $S^{n-1}=\left \{ x\in {\mathbb{R}^{n}:}\text{ }|x|=1\right \} $ denote
the unit sphere on ${\mathbb{R}^{n}}$ $(n\geq 2)$ equipped with the
normalized Lebesgue measure $d\sigma \left( x^{\prime }\right) $, where $%
x^{\prime }$ denotes the unit vector in the direction of $x$.

To study the existence and regularity results for an elliptic differential
operator, i.e.%
\begin{equation*}
D=\dsum \limits_{i,j=1}^{n}a_{i,j}\frac{\partial ^{2}}{\partial
x_{i}\partial x_{j}}
\end{equation*}%
with constant coefficients $\left \{ a_{i,j}\right \} $, among some other
estimates, one needs to study the singular integral operator $\overline{T}$
with a convolution kernel $K$ (see \cite{Calderon and Zygmund1} or \cite%
{Calderon and Zygmund2}) satisfying

$(a)$ $K\left( tx_{1},\ldots ,tx_{n}\right) =t^{-n}K\left( x\right) $,$%
\qquad $for any $t>0$;

$(b)$ $K\in C^{\infty }\left( {\mathbb{R}^{n}}\setminus \{0\} \right) $;

$(c)\int \limits_{S^{n-1}}K(x^{\prime })d\sigma (x^{\prime })=0$.

Similarly, for the heat operator%
\begin{equation*}
D=\frac{\partial }{\partial x_{1}}-\dsum \limits_{j=2}^{n}\frac{\partial ^{2}%
}{\partial x_{j}^{2}},
\end{equation*}%
the corresponding singular integral operator $\overline{T}$ have a kernel $K$
satisfying

$(a^{\prime })$ $K\left( t^{2}x_{1},\ldots ,tx_{n}\right) =t^{-n-1}K\left(
x\right) $,$\qquad $for any $t>0$;

$(b^{\prime })$ $K\in C^{\infty }\left( {\mathbb{R}^{n}}\setminus \{0\}
\right) $;

$(c^{\prime })\int \limits_{S^{n-1}}K(x^{\prime })\left( 2x_{1}^{\prime
2}+x_{2}^{\prime 2}+\cdots +x_{n}^{\prime 2}\right) d\sigma (x^{\prime })=0$.

To study the regularity results for a more general parabolic differential
operator with constant coefficients, in 1966, Fabes and Rivi\'{e}re \cite%
{Fabes and Riviere} introduced the following parabolic singular integral
operator%
\begin{equation*}
\overline{T}^{P}f(x)=p.v.\int \limits_{{\mathbb{R}^{n}}}K(y)f(x-y)\,dy
\end{equation*}%
with $K$ satisfying

$(i)$ $K\left( t^{\alpha _{1}}x_{1},\ldots ,t^{\alpha _{n}}x_{n}\right)
=t^{-\alpha }K\left( x_{1},\ldots ,x_{n}\right) $, $t>0$, $x\neq 0$, $\alpha
=\dsum \limits_{i=1}^{n}\alpha _{i}$;

$(ii)$ $K\in C^{\infty }\left( {\mathbb{R}^{n}}\setminus \{0\} \right) $;

$(iii)\int \limits_{S^{n-1}}K(x^{\prime })J\left( x^{\prime }\right) d\sigma
(x^{\prime })=0$, where $\alpha _{i}\geq 1$ $\left( i=1,\ldots ,n\right) $
and $J\left( x^{\prime }\right) =\alpha _{1}x_{1}^{\prime 2}+\ldots +\alpha
_{n}x_{n}^{\prime 2}$ is shown as follows.

Let $\rho \in \left( 0,\infty \right) $ and $0\leq \varphi _{n-1}\leq 2\pi $%
, $0\leq \varphi _{i}\leq \pi $, $i=1,\ldots ,n-2$. For any $x\in {\mathbb{R}%
^{n}}$, set%
\begin{eqnarray*}
x_{1} &=&\rho ^{\alpha _{1}}\cos \varphi _{1}\ldots \cos \varphi _{n-2}\cos
\varphi _{n-1}, \\
x_{2} &=&\rho ^{\alpha _{2}}\cos \varphi _{1}\ldots \cos \varphi _{n-2}\cos
\varphi _{n-1}, \\
&&\vdots \\
x_{n-1} &=&\rho ^{\alpha _{n-1}}\cos \varphi _{1}\sin \varphi _{2}, \\
x_{n} &=&\rho ^{\alpha _{n}}\sin \varphi _{1}.
\end{eqnarray*}

Then $dx=\rho ^{\alpha -1}J\left( x^{\prime }\right) d\rho d\sigma
(x^{\prime })$, where $\alpha =\dsum \limits_{i=1}^{n}\alpha _{i}$, $%
x^{\prime }\in S^{n-1}$, $d\sigma $ is the element of area of $S^{n-1}$ and $%
\rho ^{\alpha -1}J$ is the Jacobian of the above transform. In \cite{Fabes
and Riviere} Fabes and Rivi\'{e}re has pointed out that $J\left( x^{\prime
}\right) $ is a $C^{\infty }$ function on $S^{n-1}$ and $1\leq J\left(
x^{\prime }\right) \leq M$, where $M$ is a constant independent of $%
x^{\prime }$. Without loss of generality, in this paper we may assume $%
\alpha _{n}\geq \alpha _{n-1}\geq \cdots \geq \alpha _{1}\geq 1$. Notice
that the above condition $\left( i\right) $ can be written as

$\left( i^{\prime }\right) $ $K\left( A_{t}x\right) =\left \vert \det \left(
A_{t}\right) \right \vert ^{-1}K\left( x\right) $, where $A_{t}=diag\left[
t^{\alpha _{1}},\ldots ,t^{\alpha _{n}}\right] =%
\begin{pmatrix}
t^{\alpha _{1}} &  & 0 \\ 
& \ddots &  \\ 
0 &  & t^{\alpha _{n}}%
\end{pmatrix}%
$ is a diagonal matrix.

Note that for each fixed $x=\left( x_{1},\ldots ,x_{n}\right) \in {\mathbb{R}%
^{n}}$, the function%
\begin{equation*}
F\left( x,\rho \right) =\dsum \limits_{i=1}^{n}\frac{x_{i}^{2}}{\rho
^{2\alpha _{i}}}
\end{equation*}%
is a strictly decreasing function of $\rho >0$. Hence, there exists an
unique $t$ such that $F\left( x,t\right) =1$. It has been proved in \cite%
{Fabes and Riviere} that if we set $\rho \left( 0\right) =0$ and $\rho
\left( x\right) =t$ such that $F\left( x,t\right) =1$, then $\rho $ is a
metric on ${\mathbb{R}^{n}}$, and $\left( {\mathbb{R}^{n},}\rho \right) $ is
called the mixed homogeneity space related to $\left \{ \alpha _{i}\right \}
_{i=1}^{n}$.

\begin{remark}
Many works have been done for parabolic singular integral operators,
including the weak type estimates and $L_{p}$ (strong $\left( p,p\right) $)
boundedness. For example, one can see references \cite{Hofmann1, Hofmann2,
Tao} for details.
\end{remark}

Let $P$ be a real $n\times n$ matrix, whose all the eigenvalues have
positive real part. Let $A_{t}=t^{P}$ $\left( t>0\right) $, and set $\gamma
=trP$. Then, there exists a quasi-distance $\rho $ associated with $P$ such
that (see \cite{Coifman-Weiss})

$\left( 1-1\right) $ $\rho \left( A_{t}x\right) =t\rho \left( x\right) $, $%
t>0$, for every $x\in {\mathbb{R}^{n}}$,

$\left( 1-2\right) $ $\rho \left( 0\right) =0$, $\rho \left( x-y\right)
=\rho \left( y-x\right) \geq 0$, and $\rho \left( x-y\right) \leq k\left(
\rho \left( x-z\right) +\rho \left( y-z\right) \right) $,

$\left( 1-3\right) $ $dx=\rho ^{\gamma -1}d\sigma \left( w\right) d\rho $,
where $\rho =\rho \left( x\right) $, $w=A_{\rho ^{-1}}x$ and $d\sigma \left(
w\right) $ is a measure on the unit ellipsoid $\left \{ w:\rho \left(
w\right) =1\right \} $.

Then, $\left \{ {\mathbb{R}^{n},\rho ,dx}\right \} $ becomes a space of
homogeneous type in the sense of Coifman-Weiss (see \cite{Coifman-Weiss})
and a homogeneous group in the sense of Folland-Stein (see \cite%
{Folland-Stein}). Moreover, we always assume that there hold the following
properties of the quasi-distance $\rho $:

$\left( 1-4\right) $ For every $x$,%
\begin{eqnarray*}
c_{1}\left \vert x\right \vert ^{\alpha _{1}} &\leq &\rho \left( x\right)
\leq c_{2}\left \vert x\right \vert ^{\alpha _{2}}\text{ if }\rho \left(
x\right) \geq 1; \\
c_{3}\left \vert x\right \vert ^{\alpha _{3}} &\leq &\rho \left( x\right)
\leq c_{4}\left \vert x\right \vert ^{\alpha _{4}}\text{ if }\rho \left(
x\right) \leq 1,
\end{eqnarray*}%
and%
\begin{equation*}
\rho \left( \theta x\right) \leq \rho \left( x\right) \text{ for }0<\theta
<1,
\end{equation*}%
with some positive constants $\alpha _{i}$ and $c_{i}$ $\left( i=1,\ldots
,4\right) $. Similar properties also hold for the quasimetric $\rho ^{\ast }$
associated with the adjoint matrix $P^{\ast }$.

The following are some important examples of the above defined matrices $P$
and distances $\rho $:

$1.$ Let $\left( Px,x\right) \geq \left( x,x\right) $ $\left( x\in {\mathbb{R%
}^{n}}\right) $. In this case, $\rho \left( x\right) $ is defined by the
unique solution of $\left \vert A_{t^{-1}}x\right \vert =1$, and $k=1$. This
is the case studied by Calder\'{o}n and Torchinsky in \cite{Calderon and
Torchinsky}.

$2.$ Let $P$ be a diagonal matrix with positive diagonal entries, and let $%
t=\rho \left( x\right) $, $x\in {\mathbb{R}^{n}}$ be the unique solution of $%
\left \vert A_{t^{-1}}x\right \vert =1$.

$a)$ When all diagonal entries are greater than or equal to $1$, Besov et
al. in \cite{Besov-Lizorkin} and Fabes and Rivi\'{e}re in \cite{Fabes and
Riviere} have studied the weak $\left( 1.1\right) $ and $L_{p}$ (strong $%
\left( p,p\right) $) estimates of the singular integral operators on this
space.

$b)$ If there are diagonal entries smaller than $1$, then $\rho $ satisfies
the above $\left( 1-1\right) -\left( 1-4\right) $ with $k\geq 1$.

It is a simple matter to check that ${\rho }\left( x-y\right) $ defines a
distance between any two points $x,y\in {\mathbb{R}^{n}}$. Thus ${\mathbb{R}%
^{n}}$, endowed with the metric $\rho $, defines a homogeneous metric space 
\cite{Besov-Lizorkin, Fabes and Riviere}. Denote by $E\left( x,r\right) $
the ellipsoid with center at $x$ and radius $r$, more precisely, $E\left(
x,r\right) =\left \{ y\in {\mathbb{R}^{n}:\rho }\left( x-y\right)
<r\right
\} $. For $k>0$, we denote $kE\left( x,r\right) =\left \{ y\in {%
\mathbb{R}^{n}:\rho }\left( x-y\right) <kr\right \} $. Moreover, by the
property of $\rho $ and the polar coordinates transform above, we have%
\begin{equation*}
\left \vert E\left( x,r\right) \right \vert =\dint \limits_{{\rho }\left(
x-y\right) <r}dy=\upsilon _{\rho }r^{\alpha _{1}+\cdots +\alpha
_{n}}=\upsilon _{\rho }r^{\gamma },
\end{equation*}%
where $|E(x,r)|$ stands for the Lebesgue measure of $E(x,r)$ and $\upsilon
_{\rho }$ is the volume of the unit ellipsoid on ${\mathbb{R}^{n}}$. By $%
E^{C}(x,r)={\mathbb{R}^{n}}\setminus $ $E\left( x,r\right) $, we denote the
complement of $E\left( x,r\right) $. Moreover, in the standard parabolic
case $P_{0}=diag\left[ 1,\ldots ,1,2\right] $ we have%
\begin{equation*}
\rho \left( x\right) =\sqrt{\frac{\left \vert x^{\prime }\right \vert ^{2}+%
\sqrt{\left \vert x^{\prime }\right \vert ^{4}+x_{n}^{2}}}{2},}\qquad
x=\left( x^{\prime },x_{n}\right) .
\end{equation*}

Note that we deal not exactly with the parabolic metric, but with a general
anisotropic metric $\rho $ of generalized homogeneity, the parabolic metric
being its particular case, but we keep the term parabolic in the title and
text of the paper, the above existing tradition, see for instance \cite%
{Calderon and Torchinsky}.

Suppose that $\Omega \left( x\right) $ is a real-valued and measurable
function defined on ${\mathbb{R}^{n}}$. Suppose that $S^{n-1}$ is the unit
sphere on ${\mathbb{R}^{n}}$ $(n\geq 2)$ equipped with the normalized
Lebesgue surface measure $d\sigma $.

Let $\Omega \in L_{s}(S^{n-1})$ with $1<s\leq \infty $ be homogeneous of
degree zero with respect to $A_{t}$ ($\Omega \left( x\right) $ is $A_{t}$%
-homogeneous of degree zero). We define $s^{\prime }=\frac{s}{s-1}$ for any $%
s>1$. Suppose that $T_{\Omega }^{P}$ represents a parabolic linear or a
parabolic sublinear operator, which satisfies that for any $f\in L_{1}({%
\mathbb{R}^{n}})$ with compact support and $x\notin suppf$ 
\begin{equation}
|T_{\Omega }^{P}f(x)|\leq c_{0}\int \limits_{{\mathbb{R}^{n}}}\frac{|\Omega
(x-y)|}{{\rho }\left( x-y\right) ^{\gamma }}\,|f(y)|\,dy,  \label{e1}
\end{equation}%
where $c_{0}$ is independent of $f$ and $x$.

We point out that the condition (\ref{e1}) in the case $\Omega \equiv 1$ and 
$P=I$ was first introduced by Soria and Weiss in \cite{SW} . The condition (%
\ref{e1}) is satisfied by many interesting operators in harmonic analysis,
such as the parabolic Calder\'{o}n--Zygmund operators, parabolic Carleson's
maximal operator, parabolic Hardy--Littlewood maximal operator, parabolic C.
Fefferman's singular multipliers, parabolic R. Fefferman's singular
integrals, parabolic Ricci--Stein's oscillatory singular integrals,
parabolic the Bochner--Riesz means and so on (see \cite{LLY}, \cite{SW} for
details).

\bigskip Let $\Omega \in L_{s}(S^{n-1})$ with $1<s\leq \infty $ be
homogeneous of degree zero with respect to $A_{t}$ ($\Omega \left( x\right) $
is $A_{t}$-homogeneous of degree zero), that is,%
\begin{equation*}
\Omega (A_{t}x)=\Omega (x),
\end{equation*}
for any$~~t>0,$ $x\in {\mathbb{R}^{n}}$ and satisfies the
cancellation(vanishing) condition 
\begin{equation*}
\dint \limits_{S^{n-1}}\Omega (x^{\prime })J\left( x^{\prime }\right)
d\sigma (x^{\prime })=0,
\end{equation*}%
where $x^{\prime }=\frac{x}{|x|}$ for any $x\neq 0$.

Let $f\in L^{loc}\left( {\mathbb{R}^{n}}\right) $. The parabolic homogeneous
singular integral operator $\overline{T}_{\Omega }^{P}$ and the parabolic
maximal operator $M_{\Omega }^{P}$ by with rough kernels are defined by

\begin{equation}
\overline{T}_{\Omega }^{P}f(x)=p.v.\int \limits_{{\mathbb{R}^{n}}}\frac{%
\Omega (x-y)}{{\rho }\left( x-y\right) ^{\gamma }}f(y)dy,  \label{1*}
\end{equation}%
\begin{equation*}
M_{\Omega }^{P}f(x)=\sup_{t>0}|E(x,t)|^{-1}\int \limits_{E(x,t)}\left \vert
\Omega \left( x-y\right) \right \vert |f(y)|dy,
\end{equation*}%
satisfy condition (\ref{e1}).

It is obvious that when $\Omega \equiv 1$, $\overline{T}_{\Omega }^{P}\equiv 
\overline{T}^{P}$ and $M_{\Omega }^{P}\equiv M^{P}$ are the parabolic
singular operator and the parabolic maximal operator, respectively. If $P=I$%
, then $M_{\Omega }^{I}\equiv M_{\Omega }$ is the Hardy-Littlewood maximal
operator with rough kernel, and $\overline{T}_{\Omega }^{I}\equiv \overline{T%
}_{\Omega }$ is the homogeneous singular integral operator. It is well known
that the parabolic maximal and singular operators play an important role in
harmonic analysis (see \cite{Calderon and Torchinsky, Folland-Stein,
Palagachev, Stein93}). In particular, the boundedness of $\overline{T}%
_{\Omega }^{P}$ on Lebesgue spaces has been obtained.

\begin{theorem}
\label{Teo-Guliyev}Suppose that $\Omega \in L_{s}(S^{n-1})$, $1<s\leq \infty 
$, is $A_{t}$-homogeneous of degree zero \textit{has mean value zero on }$%
S^{n-1}$\textit{.} If $s^{\prime }\leq p$ or $p<s$, then the operator $%
\overline{T}_{\Omega }^{P}$ is bounded on $L_{p}({\mathbb{R}^{n}})$. \textit{%
Also,} the operator $\overline{T}_{\Omega }^{P}$ is bounded from $L_{1}({%
\mathbb{R}^{n}})$ to $WL_{1}({\mathbb{R}^{n}})$. Moreover, we have for $p>1$

\begin{equation*}
\left \Vert \overline{T}_{\Omega }^{P}f\right \Vert _{L_{p}}\leq C\left
\Vert f\right \Vert _{L_{p}},
\end{equation*}%
and for $p=1$%
\begin{equation*}
\left \Vert \overline{T}_{\Omega }^{P}f\right \Vert _{WL_{1}}\leq C\left
\Vert f\right \Vert _{L_{1}}.
\end{equation*}
\end{theorem}

\begin{corollary}
Under the assumptions of Theorem \ref{Teo-Guliyev}, the operator $M_{\Omega
}^{P}$ is bounded on $L_{p}({\mathbb{R}^{n}})$. \textit{Also,} the operator $%
M_{\Omega }^{P}$ is bounded from $L_{1}({\mathbb{R}^{n}})$ to $WL_{1}({%
\mathbb{R}^{n}})$. Moreover, we have for $p>1$

\begin{equation*}
\left \Vert M_{\Omega }^{P}f\right \Vert _{L_{p}}\leq C\left \Vert f\right
\Vert _{L_{p}},
\end{equation*}%
and for $p=1$%
\begin{equation*}
\left \Vert M_{\Omega }^{P}f\right \Vert _{WL_{1}}\leq C\left \Vert f\right
\Vert _{L_{1}}.
\end{equation*}
\end{corollary}

\begin{proof}
It suffices to refer to the known fact that%
\begin{equation*}
M_{\Omega }^{P}f\left( x\right) \leq C_{\gamma }\overline{T}_{\Omega
}^{P}f\left( x\right) ,\qquad C_{\gamma }=\left \vert E\left( 0,1\right)
\right \vert .
\end{equation*}
\end{proof}

Note that in the isotropic case $P=I$ Theorem \ref{Teo-Guliyev} has been
proved in \cite{Muckenhoupt}.

Let $b$ be a locally integrable function on ${\mathbb{R}^{n}}$, then we
define commutators generated by parabolic maximal and singular integral
operators by with rough kernels and $b$ as follows, respectively.%
\begin{equation*}
M_{\Omega ,b}^{P}\left( f\right) (x)=\sup_{t>0}|E(x,t)|^{-1}\int
\limits_{E(x,t)}\left \vert b\left( x\right) -b\left( y\right) \right \vert
\left \vert \Omega \left( x-y\right) \right \vert |f(y)|dy,
\end{equation*}%
\begin{equation}
\lbrack b,\overline{T}_{\Omega }^{P}]f(x)\equiv b(x)\overline{T}_{\Omega
}^{P}f(x)-\overline{T}_{\Omega }^{P}(bf)(x)=p.v.\int \limits_{{\mathbb{R}^{n}%
}}[b(x)-b(y)]\frac{\Omega (x-y)}{{\rho }\left( x-y\right) ^{\gamma }}f(y)dy.
\label{2*}
\end{equation}

If we take $\alpha _{1}=\cdots \alpha _{n}=1$ and $P=I$, then obviously $%
\rho \left( x\right) =\left \vert x\right \vert =\left( \dsum
\limits_{i=1}^{n}x_{i}^{2}\right) ^{\frac{1}{2}}$, $\gamma =n$, $\left( {%
\mathbb{R}^{n},\rho }\right) =$ $\left( {\mathbb{R}^{n},}\left \vert \cdot
\right \vert \right) $, $E_{I}(x,r)=B\left( x,r\right) $, $A_{t}=tI$ and $%
J\left( x^{\prime }\right) \equiv 1$. In this case, $\overline{T}_{\Omega
}^{P}$ defined as in (\ref{1*}) is the classical singular integral operator
with rough kernel of convolution type whose boundedness in various function
spaces has been well-studied by many authors (see \cite{BGGS, Lu1}, \cite%
{Gurbuz1, Gurbuz2}, \cite{LuDingY}, and so on). And also, in this case, $[b,%
\overline{T}_{\Omega }^{P}]$ defined as in (\ref{2*}) is the classical
commutator of singular integral operator with rough kernel of convolution
type whose boundedness in various function spaces has also been well-studied
by many authors (see \cite{BGGS, Lu1}, \cite{Gurbuz1, Gurbuz2}, \cite%
{LuDingY}, and so on).

The classical Morrey spaces $M_{p,\lambda }$ have been introduced by Morrey
in \cite{Morrey} to study the local behavior of solutions of second order
elliptic partial differential equations(PDEs). In recent years there has
been an explosion of interest in the study of the boundedness of operators
on Morrey-type spaces. It has been obtained that many properties of
solutions to PDEs are concerned with the boundedness of some operators on
Morrey-type spaces. In fact, better inclusion between Morrey and H\"{o}lder
spaces allows to obtain higher regularity of the solutions to different
elliptic and parabolic boundary problems.

Morrey has stated that many properties of solutions to PDEs can be
attributed to the boundedness of some operators on Morrey spaces. For the
boundedness of the Hardy--Littlewood maximal operator, the fractional
integral operator and the Calder\'{o}n--Zygmund singular integral operator
on these spaces, we refer the readers to \cite{Adams, ChFra, Peetre}. For
the properties and applications of classical Morrey spaces, see \cite%
{ChFraL1, ChFraL2, FazRag2, FazPalRag} and references therein. The
generalized Morrey spaces $M_{p,\varphi }$ are obtained by replacing $%
r^{\lambda }$ with a function $\varphi \left( r\right) $ in the definition
of the Morrey space. During the last decades various classical operators,
such as maximal, singular and potential operators have been widely
investigated in classical and generalized Morrey spaces.

We define the parabolic Morrey spaces $M_{p,\lambda ,P}\left( {\mathbb{R}^{n}%
}\right) $ via the norm

\begin{equation*}
\left \Vert f\right \Vert _{M_{p,\lambda ,P}}=\sup \limits_{x\in {\mathbb{R}%
^{n}},r>0}\,r^{-\frac{\lambda }{p}}\, \Vert f\Vert _{L_{p}(E(x,r))}<\infty ,
\end{equation*}%
where $f\in L_{p}^{loc}({\mathbb{R}^{n}})$, $0\leq \lambda \leq \gamma $ and 
$1\leq p\leq \infty $.

Note that $M_{p,0,P}=L_{p}({\mathbb{R}^{n}})$ and $M_{p,\gamma ,P}=L_{\infty
}({\mathbb{R}^{n}})$. If $\lambda <0$ or $\lambda >\gamma $, then $%
M_{p,\lambda }={\Theta }$, where $\Theta $ is the set of all functions
equivalent to $0$ on ${\mathbb{R}^{n}}$.

We also denote by $WM_{p,\lambda ,P}\equiv WM_{p,\lambda ,P}({\mathbb{R}^{n}}%
)$ the weak parabolic Morrey space of all functions $f\in WL_{p}^{loc}({%
\mathbb{R}^{n}})$ for which 
\begin{equation*}
\left \Vert f\right \Vert _{WM_{p,\lambda ,P}}\equiv \left \Vert f\right
\Vert _{WM_{p,\lambda ,P}({\mathbb{R}^{n}})}=\sup_{x\in {\mathbb{R}^{n}}%
,r>0}r^{-\frac{\lambda }{p}}\Vert f\Vert _{WL_{p}(E(x,r))}<\infty ,
\end{equation*}%
where $WL_{p}(E(x,r))$ denotes the weak $L_{p}$-space of measurable
functions $f$ for which 
\begin{equation*}
\begin{split}
\Vert f\Vert _{WL_{p}(E(x,r))}& \equiv \Vert f\chi _{_{E(x,r)}}\Vert
_{WL_{p}({\mathbb{R}^{n}})} \\
& =\sup_{t>0}t\left \vert \left \{ y\in E(x,r):\,|f(y)|>t\right \} \right
\vert ^{1/{p}} \\
& =\sup_{0<t\leq |E(x,r)|}t^{1/{p}}\left( f\chi _{_{E(x,r)}}\right) ^{\ast
}(t)<\infty ,
\end{split}%
\end{equation*}%
where $g^{\ast }$ denotes the non-increasing rearrangement of a function $g$.

Note that $WL_{p}({\mathbb{R}^{n}})=WM_{p,0,P}\left( {\mathbb{R}^{n}}\right) 
$,%
\begin{equation*}
M_{p,\lambda ,P}({\mathbb{R}^{n}})\subset WM_{p,\lambda ,P}({\mathbb{R}^{n}})%
\text{ and }\left \Vert f\right \Vert _{WM_{p,\lambda ,P}}\leq \left \Vert
f\right \Vert _{M_{p,\lambda ,P}}.
\end{equation*}

If $P=I$, then $M_{p,\lambda ,I}({\mathbb{R}^{n}})\equiv M_{p,\lambda }({%
\mathbb{R}^{n}})$ is the classical Morrey space.

It is known that the parabolic maximal operator $M^{P}$ is also bounded on $%
M_{p,\lambda ,P}$ for all $1<p<\infty $ and $0<\lambda <\gamma $ (see, e.g. 
\cite{Meskhi}), whose isotropic counterpart has been proved by Chiarenza and
Frasca \cite{ChFra}.

In this paper, we prove the boundedness of the parabolic sublinear operators
with rough kernel $T_{\Omega }^{P}$ satisfying condition (\ref{e1})
generated by parabolic Calder\'{o}n-Zygmund operators with rough kernel from
one parabolic generalized local Morrey space $LM_{p,\varphi
_{1},P}^{\{x_{0}\}}$ to another one $LM_{p,\varphi _{2},P}^{\{x_{0}\}}$, $%
1<p<\infty $, and from the space $LM_{1,\varphi _{1},P}^{\{x_{0}\}}$ to the
weak space $WLM_{1,\varphi _{2},P}^{\{x_{0}\}}$. In the case of $b\in
LC_{p_{2},\lambda ,P}^{\left \{ x_{0}\right \} }$ (parabolic local Campanato
space) and $[b,T_{\Omega }^{P}]$ is a sublinear operator, we find the
sufficient conditions on the pair $(\varphi _{1},\varphi _{2})$ which
ensures the boundedness of the commutator operators $[b,T_{\Omega }^{P}]$
from $LM_{p_{1},\varphi _{1},P}^{\{x_{0}\}}$ to $LM_{p,\varphi
_{2},P}^{\{x_{0}\}}$, $1<p<\infty $, $\frac{1}{p}=\frac{1}{p_{1}}+\frac{1}{%
p_{2}}$ and $0\leq \lambda <\frac{1}{\gamma }$.

By $A\lesssim B$ we mean that $A\leq CB$ with some positive constant $C$
independent of appropriate quantities. If $A\lesssim B$ and $B\lesssim A$,
we write $A\approx B$ and say that $A$ and $B$ are equivalent.

\section{parabolic generalized local Morrey spaces}

Let us define the parabolic generalized Morrey spaces as follows.

\begin{definition}
\label{Definition1}\textbf{(parabolic generalized Morrey space) }Let $%
\varphi (x,r)$ be a positive measurable function on ${\mathbb{R}^{n}}\times
(0,\infty )$ and $1\leq p<\infty $. We denote by $M_{p,\varphi ,P}\equiv
M_{p,\varphi ,P}({\mathbb{R}^{n}})$ the parabolic generalized Morrey space,
the space of all functions $f\in L_{p}^{loc}({\mathbb{R}^{n}})$ with finite
quasinorm 
\begin{equation*}
\Vert f\Vert _{M_{p,\varphi ,P}}=\sup \limits_{x\in {\mathbb{R}^{n}}%
,r>0}\varphi (x,r)^{-1}\,|E(x,r)|^{-\frac{1}{p}}\, \Vert f\Vert
_{L_{p}(E(x,r))}<\infty .
\end{equation*}%
Also by $WM_{p,\varphi ,P}\equiv WM_{p,\varphi ,P}({\mathbb{R}^{n}})$ we
denote the weak parabolic generalized Morrey space of all functions $f\in
WL_{p}^{loc}({\mathbb{R}^{n}})$ for which 
\begin{equation*}
\Vert f\Vert _{WM_{p,\varphi ,P}}=\sup \limits_{x\in {\mathbb{R}^{n}}%
,r>0}\varphi (x,r)^{-1}\,|E(x,r)|^{-\frac{1}{p}}\, \Vert f\Vert
_{WL_{p}(E(x,r))}<\infty .
\end{equation*}
\end{definition}

According to this definition, we recover the parabolic Morrey space $%
M_{p,\lambda ,P}$ and the weak parabolic Morrey space $WM_{p,\lambda ,P}$
under the choice $\varphi (x,r)=r^{\frac{\lambda -\gamma }{p}}$: 
\begin{equation*}
M_{p,\lambda ,P}=M_{p,\varphi ,P}\mid _{\varphi (x,r)=r^{\frac{\lambda
-\gamma }{p}}},~~~~~~~~WM_{p,\lambda ,P}=WM_{p,\varphi ,P}\mid _{\varphi
(x,r)=r^{\frac{\lambda -\gamma }{p}}}.
\end{equation*}

Inspired by the above Definition \ref{Definition1}, \cite{BGGS} and the
Ph.D. thesis of Gurbuz \cite{Gurbuz1}, we introduce the parabolic
generalized local Morrey spaces $LM_{p,\varphi ,P}^{\{x_{0}\}}$ by the
following definition.

\begin{definition}
\label{Definition2}\textbf{(parabolic generalized local Morrey space) }Let $%
\varphi (x,r)$ be a positive measurable function on ${\mathbb{R}^{n}}\times
(0,\infty )$ and $1\leq p<\infty $. For any fixed $x_{0}\in {\mathbb{R}^{n}}$
we denote by $LM_{p,\varphi ,P}^{\{x_{0}\}}\equiv LM_{p,\varphi
,P}^{\{x_{0}\}}({\mathbb{R}^{n}})$ the parabolic generalized local Morrey
space, the space of all functions $f\in L_{p}^{loc}({\mathbb{R}^{n}})$ with
finite quasinorm 
\begin{equation*}
\Vert f\Vert _{LM_{p,\varphi ,P}^{\{x_{0}\}}}=\sup \limits_{r>0}\varphi
(x_{0},r)^{-1}\,|E(x_{0},r)|^{-\frac{1}{p}}\, \Vert f\Vert
_{L_{p}(E(x_{0},r))}<\infty .
\end{equation*}%
Also by $WLM_{p,\varphi ,P}^{\{x_{0}\}}\equiv WLM_{p,\varphi ,P}^{\{x_{0}\}}(%
{\mathbb{R}^{n}})$ we denote the weak parabolic generalized local Morrey
space of all functions $f\in WL_{p}^{loc}({\mathbb{R}^{n}})$ for which 
\begin{equation*}
\Vert f\Vert _{WLM_{p,\varphi ,P}^{\{x_{0}\}}}=\sup \limits_{r>0}\varphi
(x_{0},r)^{-1}\,|E(x_{0},r)|^{-\frac{1}{p}}\, \Vert f\Vert
_{WL_{p}(E(x_{0},r))}<\infty .
\end{equation*}
\end{definition}

According to this definition, we recover the local parabolic Morrey space $%
LM_{p,\lambda ,P}^{\{x_{0}\}}$ and weak local parabolic Morrey space $%
WLM_{p,\lambda ,P}^{\{x_{0}\}}$ under the choice $\varphi (x_{0},r)=r^{\frac{%
\lambda -\gamma }{p}}$:%
\begin{equation*}
LM_{p,\lambda ,P}^{\{x_{0}\}}=LM_{p,\varphi ,P}^{\{x_{0}\}}\mid _{\varphi
(x_{0},r)=r^{\frac{\lambda -\gamma }{p}}},~~~~~~WLM_{p,\lambda
,P}^{\{x_{0}\}}=WLM_{p,\varphi ,P}^{\{x_{0}\}}\mid _{\varphi (x_{0},r)=r^{%
\frac{\lambda -\gamma }{p}}}.
\end{equation*}

Furthermore, we have the following embeddings:%
\begin{equation*}
M_{p,\varphi ,P}\subset LM_{p,\varphi ,P}^{\{x_{0}\}},\qquad \Vert f\Vert
_{LM_{p,\varphi ,P}^{\{x_{0}\}}}\leq \Vert f\Vert _{M_{p,\varphi ,P}},
\end{equation*}%
\begin{equation*}
WM_{p,\varphi ,P}\subset WLM_{p,\varphi ,P}^{\{x_{0}\}},\qquad \Vert f\Vert
_{WLM_{p,\varphi ,P}^{\{x_{0}\}}}\leq \Vert f\Vert _{WM_{p,\varphi ,P}}.
\end{equation*}

In \cite{DingYZ} the following statement has been proved for parabolic
singular operators with rough kernel $\overline{T}_{\Omega }^{P}$,
containing the result in \cite{Miz, Nakai1, Nakai2}.

\begin{theorem}
Suppose that $\Omega \in L_{s}(S^{n-1})$, $1<s\leq \infty $, is $A_{t}$%
-homogeneous of degree zero and has mean value zero on $S^{n-1}$. Let $1\leq
s^{\prime }<p<\infty \left( s^{\prime }=\frac{s}{s-1}\right) $ and $\varphi
(x,r)$ satisfies conditions 
\begin{equation}
c^{-1}\varphi (x,r)\leq \varphi (x,t)\leq c\, \varphi (x,r)  \label{e32}
\end{equation}%
whenever $r\leq t\leq 2r$, where $c~(\geq 1)$ does not depend on $t$, $r$, $%
x\in {\mathbb{R}^{n}}$ and%
\begin{equation}
\int \limits_{r}^{\infty }\varphi (x,t)^{p}\frac{dt}{t}\leq C\, \varphi
(x,r)^{p},  \label{e33x}
\end{equation}%
where $C$ does not depend on $x$ and $r$. Then the operator $\overline{T}%
_{\Omega }^{P}$ is bounded on $M_{p,\varphi ,P}$.
\end{theorem}

The results of \cite{Miz, Nakai1, Nakai2} imply the following statement.

\begin{theorem}
Let $1\leq p<\infty $ and $\varphi (x,t)$ satisfies conditions (\ref{e32})
and (\ref{e33x}). Then the operators $M^{P}$ and $\overline{T}^{P}$ are
bounded on $M_{p,\varphi ,P}$ for $p>1$ and from $M_{1,\varphi ,P}$ to $%
WM_{1,\varphi ,P}$ and for $p=1$.
\end{theorem}

The following statement, containing the results obtained in \cite{Miz}, \cite%
{Nakai1}, \cite{Nakai2} has been proved in \cite{GulDoc, GulJIA} (see also 
\cite{BurGulHus1}-\cite{BurGogGulMus2}, \cite{GulDoc}-\cite{GULAKShIEOT2012}%
).

\begin{theorem}
\label{teo4}Let $1\leq p<\infty $ and the pair $(\varphi _{1},\varphi _{2})$
satisfies the condition 
\begin{equation}
\int \limits_{r}^{\infty }\varphi _{1}(x,t)\frac{dt}{t}\leq C\, \varphi
_{2}(x,r),  \label{e35}
\end{equation}%
where $C$ does not depend on $x$ and $r$. Then the operator $\overline{T}%
^{P} $ is bounded from $M_{p,\varphi _{1},P}$ to $M_{p,\varphi _{2},P}$ for $%
p>1$ and from $M_{1,\varphi _{1},P}$ to $WM_{1,\varphi _{2},P}$ for $p=1$.
\end{theorem}

Finally, inspired by the Definition \ref{Definition2}, \cite{BGGS} and the
Ph.D. thesis of Gurbuz \cite{Gurbuz1} in this paper we consider the
boundedness of parabolic sublinear operators with rough kernel on the
parabolic generalized local Morrey spaces and give the parabolic local
Campanato space estimates for their commutators.

\section{Parabolic sublinear operators with rough kernel generated by
parabolic Calder\'{o}n-Zygmund operators on the spaces $LM_{p,\protect%
\varphi ,P}^{\{x_{0}\}}$}

In this section, we will prove the boundedness of the operator $T_{\Omega
}^{P}$ on the parabolic generalized local Morrey spaces $LM_{p,\varphi
,P}^{\{x_{0}\}}$ by using the following statement on the boundedness of the
weighted Hardy operator 
\begin{equation*}
H_{\omega }g(t):=\int \limits_{t}^{\infty }g(s)\omega (s)ds,\qquad
0<t<\infty ,
\end{equation*}%
where $\omega $ is a fixed non-negative function and measurable on $%
(0,\infty )$.

\begin{theorem}
\label{teo5}\cite{BGGS, Gurbuz1} Let $v_{1}$, $v_{2}$ and $\omega $ be
positive almost everywhere and measurable functions on $(0,\infty )$. The
inequality 
\begin{equation}
\limfunc{esssup}\limits_{t>0}v_{2}(t)H_{\omega }g(t)\leq C\limfunc{esssup}%
\limits_{t>0}v_{1}(t)g(t)  \label{e31}
\end{equation}%
holds for some $C>0$ for all non-negative and non-decreasing functions $g$
on $(0,\infty )$ if and only if 
\begin{equation}
B:=\sup \limits_{t>0}v_{2}(t)\int \limits_{t}^{\infty }\frac{\omega (s)ds}{%
\limfunc{esssup}\limits_{s<\tau <\infty }v_{1}(\tau )}<\infty .  \label{3*}
\end{equation}

Moreover, the value $C=B$ is the best constant for (\ref{e31}).
\end{theorem}

We first prove the following Theorem \ref{teo2*} (our main result).

\begin{theorem}
\label{teo2*}Let $x_{0}\in {\mathbb{R}^{n}}$, $1\leq p<\infty $ and $\Omega
\in L_{s}(S^{n-1})$, $1<s\leq \infty $, be $A_{t}$-homogeneous of degree
zero. Let $T_{\Omega }^{P}$ be a parabolic sublinear operator satisfying
condition (\ref{e1}), bounded on $L_{p}({\mathbb{R}^{n}})$ for $p>1$, and
bounded from $L_{1}({\mathbb{R}^{n}})$ to $WL_{1}({\mathbb{R}^{n}})$.

\textit{If }$p>1$\textit{\ and }$s^{\prime }\leq p$\textit{, then the
inequality }%
\begin{equation*}
\left \Vert T_{\Omega }^{P}f\right \Vert _{L_{p}\left( E\left(
x_{0},r\right) \right) }\lesssim r^{\frac{\gamma }{p}}\int
\limits_{2kr}^{\infty }t^{-\frac{\gamma }{p}-1}\left \Vert f\right \Vert
_{L_{p}\left( E\left( x_{0},t\right) \right) }dt
\end{equation*}%
holds for any ellipsoid $E\left( x_{0},r\right) $ and for all $f\in
L_{p}^{loc}\left( {\mathbb{R}^{n}}\right) $.

If $p>1$ and $p<s$, then the inequality%
\begin{equation*}
\left \Vert T_{\Omega }^{P}f\right \Vert _{L_{p}\left( E\left(
x_{0},r\right) \right) }\lesssim r^{\frac{\gamma }{p}-\frac{\gamma }{s}}\int
\limits_{2kr}^{\infty }t^{\frac{\gamma }{s}-\frac{\gamma }{p}-1}\left \Vert
f\right \Vert _{L_{p}\left( E\left( x_{0},t\right) \right) }dt
\end{equation*}%
holds for any ellipsoid $E\left( x_{0},r\right) $ and for all $f\in
L_{p}^{loc}\left( {\mathbb{R}^{n}}\right) $.

Moreover, for $s>1$ the inequality%
\begin{equation}
\left \Vert T_{\Omega }^{P}f\right \Vert _{WL_{q}\left( E\left(
x_{0},r\right) \right) }\lesssim r^{\gamma }\int \limits_{2kr}^{\infty
}t^{-\gamma -1}\left \Vert f\right \Vert _{L_{1}\left( E\left(
x_{0},t\right) \right) }dt  \label{38}
\end{equation}%
holds for any ellipsoid $E\left( x_{0},r\right) $ and for all $f\in
L_{1}^{loc}\left( {\mathbb{R}^{n}}\right) $.
\end{theorem}

\begin{proof}
Let $1<p<\infty $ and $s^{\prime }\leq p$. Set $E=E\left( x_{0},r\right) $
for the parabolic ball (ellipsoid) centered at $x_{0}$ and of radius $r$ and 
$2kE=E\left( x_{0},2kr\right) $. We represent $f$ as%
\begin{equation*}
f=f_{1}+f_{2},\qquad \text{\ }f_{1}\left( y\right) =f\left( y\right) \chi
_{2kE}\left( y\right) ,\qquad \text{\ }f_{2}\left( y\right) =f\left(
y\right) \chi _{\left( 2kE\right) ^{C}}\left( y\right) ,\qquad r>0
\end{equation*}%
and have%
\begin{equation*}
\left \Vert T_{\Omega }^{P}f\right \Vert _{L_{p}\left( E\right) }\leq \left
\Vert T_{\Omega }^{P}f_{1}\right \Vert _{L_{p}\left( E\right) }+\left \Vert
T_{\Omega }^{P}f_{2}\right \Vert _{L_{p}\left( E\right) }.
\end{equation*}

Since $f_{1}\in L_{p}\left( \mathbb{R}^{n}\right) $, $T_{\Omega
}^{P}f_{1}\in L_{p}\left( \mathbb{R}^{n}\right) $ and from the boundedness
of $T_{\Omega }^{P}$ on $L_{p}({\mathbb{R}^{n}})$ (see Theorem \ref%
{Teo-Guliyev}) it follows that:%
\begin{equation*}
\left \Vert T_{\Omega }^{P}f_{1}\right \Vert _{L_{p}\left( E\right) }\leq
\left \Vert T_{\Omega }^{P}f_{1}\right \Vert _{L_{p}\left( 
\mathbb{R}
^{n}\right) }\leq C\left \Vert f_{1}\right \Vert _{L_{p}\left( 
\mathbb{R}
^{n}\right) }=C\left \Vert f\right \Vert _{L_{p}\left( 2kE\right) },
\end{equation*}%
where constant $C>0$ is independent of $f$.

It is clear that $x\in E$, $y\in \left( 2kE\right) ^{C}$ implies $\frac{1}{2k%
}\rho \left( x_{0}-y\right) \leq \rho \left( x-y\right) \leq \frac{3k}{2}%
\rho \left( x_{0}-y\right) $. We get%
\begin{equation*}
\left \vert T_{\Omega }^{P}f_{2}\left( x\right) \right \vert \leq 2^{\gamma
}c_{1}\int \limits_{\left( 2kE\right) ^{C}}\frac{\left \vert f\left(
y\right) \right \vert \left \vert \Omega \left( x-y\right) \right \vert }{%
\rho \left( x_{0}-y\right) ^{\gamma }}dy.
\end{equation*}

By the Fubini's theorem, we have%
\begin{eqnarray*}
\int \limits_{\left( 2kE\right) ^{C}}\frac{\left \vert f\left( y\right)
\right \vert \left \vert \Omega \left( x-y\right) \right \vert }{\rho \left(
x_{0}-y\right) ^{\gamma }}dy &\approx &\int \limits_{\left( 2kE\right)
^{C}}\left \vert f\left( y\right) \right \vert \left \vert \Omega \left(
x-y\right) \right \vert \int \limits_{\rho \left( x_{0}-y\right) }^{\infty }%
\frac{dt}{t^{\gamma +1}}dy \\
&\approx &\int \limits_{2kr}^{\infty }\int \limits_{2kr\leq \rho \left(
x_{0}-y\right) \leq t}\left \vert f\left( y\right) \right \vert \left \vert
\Omega \left( x-y\right) \right \vert dy\frac{dt}{t^{\gamma +1}} \\
&\lesssim &\int \limits_{2kr}^{\infty }\int \limits_{E\left( x_{0},t\right)
}\left \vert f\left( y\right) \right \vert \left \vert \Omega \left(
x-y\right) \right \vert dy\frac{dt}{t^{\gamma +1}}.
\end{eqnarray*}

Applying the H\"{o}lder's inequality, we get%
\begin{eqnarray}
&&\int \limits_{\left( 2kE\right) ^{C}}\frac{\left \vert f\left( y\right)
\right \vert \left \vert \Omega \left( x-y\right) \right \vert }{\rho \left(
x_{0}-y\right) ^{\gamma }}dy  \notag \\
&\lesssim &\int \limits_{2kr}^{\infty }\left \Vert f\right \Vert
_{L_{p}\left( E\left( x_{0},t\right) \right) }\left \Vert \Omega \left(
x-\cdot \right) \right \Vert _{L_{s}\left( E\left( x_{0},t\right) \right)
}\left \vert E\left( x_{0},t\right) \right \vert ^{1-\frac{1}{p}-\frac{1}{s}}%
\frac{dt}{t^{\gamma +1}}.  \label{e310}
\end{eqnarray}

For $x\in E\left( x_{0},t\right) $, notice that $\Omega $ is $A_{t}$%
-homogenous of degree zero and $\Omega \in L_{s}(S^{n-1})$, $s>1$. Then, we
obtain%
\begin{eqnarray}
\left( \int \limits_{E\left( x_{0},t\right) }\left \vert \Omega \left(
x-y\right) \right \vert ^{s}dy\right) ^{\frac{1}{s}} &=&\left( \int
\limits_{E\left( x-x_{0},t\right) }\left \vert \Omega \left( z\right) \right
\vert ^{s}dz\right) ^{\frac{1}{s}}  \notag \\
&\leq &\left( \int \limits_{E\left( 0,t+\left \vert x-x_{0}\right \vert
\right) }\left \vert \Omega \left( z\right) \right \vert ^{s}dz\right) ^{%
\frac{1}{s}}  \notag \\
&\leq &\left( \int \limits_{E\left( 0,2t\right) }\left \vert \Omega \left(
z\right) \right \vert ^{s}dz\right) ^{\frac{1}{s}}  \notag \\
&=&\left( \int \limits_{S^{n-1}}\int \limits_{0}^{2t}\left \vert \Omega
\left( z^{\prime }\right) \right \vert ^{s}d\sigma \left( z^{\prime }\right)
r^{n-1}dr\right) ^{\frac{1}{s}}  \notag \\
&=&C\left \Vert \Omega \right \Vert _{L_{s}\left( S^{n-1}\right) }\left
\vert E\left( x_{0},2t\right) \right \vert ^{\frac{1}{s}}.  \label{e311}
\end{eqnarray}

Thus, by (\ref{e311}), it follows that:%
\begin{equation*}
\left \vert T_{\Omega }^{P}f_{2}\left( x\right) \right \vert \lesssim \int
\limits_{2kr}^{\infty }\left \Vert f\right \Vert _{L_{p}\left( E\left(
x_{0},t\right) \right) }\frac{dt}{t^{\frac{\gamma }{p}+1}}.
\end{equation*}

Moreover, for all $p\in \left[ 1,\infty \right) $ the inequality%
\begin{equation}
\left \Vert T_{\Omega }^{P}f_{2}\right \Vert _{L_{p}\left( E\right)
}\lesssim r^{\frac{\gamma }{p}}\int \limits_{2kr}^{\infty }\left \Vert
f\right \Vert _{L_{p}\left( E\left( x_{0},t\right) \right) }\frac{dt}{t^{%
\frac{\gamma }{p}+1}}  \label{e312}
\end{equation}

is valid. Thus, we obtain%
\begin{equation*}
\left \Vert T_{\Omega }^{P}f\right \Vert _{L_{p}\left( E\right) }\lesssim
\left \Vert f\right \Vert _{L_{p}\left( 2kE\right) }+r^{\frac{\gamma }{p}%
}\int \limits_{2kr}^{\infty }\left \Vert f\right \Vert _{L_{p}\left( E\left(
x_{0},t\right) \right) }\frac{dt}{t^{\frac{\gamma }{p}+1}}.
\end{equation*}

On the other hand, we have%
\begin{eqnarray}
\left \Vert f\right \Vert _{L_{p}\left( 2kE\right) } &\approx &r^{\frac{%
\gamma }{p}}\left \Vert f\right \Vert _{L_{p}\left( 2kE\right) }\int
\limits_{2kr}^{\infty }\frac{dt}{t^{\frac{\gamma }{p}+1}}  \notag \\
&\leq &r^{\frac{\gamma }{p}}\int \limits_{2kr}^{\infty }\left \Vert f\right
\Vert _{L_{p}\left( E\left( x_{0},t\right) \right) }\frac{dt}{t^{\frac{%
\gamma }{p}+1}}.  \label{e313}
\end{eqnarray}

By combining the above inequalities, we obtain%
\begin{equation*}
\left \Vert T_{\Omega }^{P}f\right \Vert _{L_{p}\left( E\right) }\lesssim r^{%
\frac{\gamma }{p}}\int \limits_{2kr}^{\infty }\left \Vert f\right \Vert
_{L_{p}\left( E\left( x_{0},t\right) \right) }\frac{dt}{t^{\frac{\gamma }{p}%
+1}}.
\end{equation*}

Let $1<p<s$. Similarly to (\ref{e311}), when $y\in B\left( x_{0},t\right) $,
it is true that%
\begin{equation}
\left( \int \limits_{E\left( x_{0},r\right) }\left \vert \Omega \left(
x-y\right) \right \vert ^{s}dy\right) ^{\frac{1}{s}}\leq C\left \Vert \Omega
\right \Vert _{L_{s}\left( S^{n-1}\right) }\left \vert E\left( x_{0},\frac{3%
}{2}t\right) \right \vert ^{\frac{1}{s}}.  \label{314}
\end{equation}

By the Fubini's theorem, the Minkowski inequality and (\ref{314}) , we get%
\begin{eqnarray*}
\left \Vert T_{\Omega }^{P}f_{2}\right \Vert _{L_{p}\left( E\right) } &\leq
&\left( \int \limits_{E}\left \vert \int \limits_{2kr}^{\infty }\int
\limits_{E\left( x_{0},t\right) }\left \vert f\left( y\right) \right \vert
\left \vert \Omega \left( x-y\right) \right \vert dy\frac{dt}{t^{\gamma +1}}%
\right \vert ^{p}dx\right) ^{\frac{1}{p}} \\
&\leq &\int \limits_{2kr}^{\infty }\int \limits_{E\left( x_{0},t\right)
}\left \vert f\left( y\right) \right \vert \left \Vert \Omega \left( \cdot
-y\right) \right \Vert _{L_{p}\left( E\right) }dy\frac{dt}{t^{\gamma +1}} \\
&\leq &\left \vert E\left( x_{0},r\right) \right \vert ^{\frac{1}{p}-\frac{1%
}{s}}\int \limits_{2kr}^{\infty }\int \limits_{E\left( x_{0},t\right) }\left
\vert f\left( y\right) \right \vert \left \Vert \Omega \left( \cdot
-y\right) \right \Vert _{L_{s}\left( E\right) }dy\frac{dt}{t^{\gamma +1}} \\
&\lesssim &r^{\frac{\gamma }{p}-\frac{\gamma }{s}}\int \limits_{2kr}^{\infty
}\left \Vert f\right \Vert _{L_{1}\left( E\left( x_{0},t\right) \right)
}\left \vert E\left( x_{0},\frac{3}{2}t\right) \right \vert ^{\frac{1}{s}}%
\frac{dt}{t^{\gamma +1}} \\
&\lesssim &r^{\frac{\gamma }{p}-\frac{\gamma }{s}}\int \limits_{2kr}^{\infty
}t^{\frac{\gamma }{s}-\frac{\gamma }{p}-1}\left \Vert f\right \Vert
_{L_{p}\left( E\left( x_{0},t\right) \right) }dt.
\end{eqnarray*}

Let $p=1<s\leq \infty $. From the weak $\left( 1,1\right) $ boundedness of $%
T_{\Omega }^{P}$ and (\ref{e313}) it follows that:%
\begin{eqnarray}
\left \Vert T_{\Omega }^{P}f_{1}\right \Vert _{WL_{1}\left( E\right) } &\leq
&\left \Vert T_{\Omega }^{P}f_{1}\right \Vert _{WL_{1}\left( 
\mathbb{R}
^{n}\right) }\lesssim \left \Vert f_{1}\right \Vert _{L_{1}\left( 
\mathbb{R}
^{n}\right) }  \notag \\
&=&\left \Vert f\right \Vert _{L_{1}\left( 2kE\right) }\lesssim r^{\gamma
}\int \limits_{2kr}^{\infty }\left \Vert f\right \Vert _{L_{1}\left( E\left(
x_{0},t\right) \right) }\frac{dt}{t^{\gamma +1}}.  \label{315}
\end{eqnarray}

Then from (\ref{e312}) and (\ref{315}) we get the inequality (\ref{38}),
which completes the proof.
\end{proof}

In the following theorem (our main result), we get the boundedness of the
operator $T_{\Omega }^{P}$ satisfying condition (\ref{e1}) on the parabolic
generalized local Morrey spaces $LM_{p,\varphi ,P}^{\{x_{0}\}}$.

\begin{theorem}
\label{teo9}Let $x_{0}\in {\mathbb{R}^{n}}$, $1\leq p<\infty $ and $\Omega
\in L_{s}(S^{n-1})$, $1<s\leq \infty $, be $A_{t}$-homogeneous of degree
zero. Let $T_{\Omega }^{P}$ be a parabolic sublinear operator satisfying
condition (\ref{e1}), bounded on $L_{p}({\mathbb{R}^{n}})$ for $p>1$, and
bounded from $L_{1}({\mathbb{R}^{n}})$ to $WL_{1}({\mathbb{R}^{n}})$. Let
also, for $s^{\prime }\leq p$, $p\neq 1$, the pair $(\varphi _{1},\varphi
_{2})$ satisfies the condition%
\begin{equation}
\int \limits_{r}^{\infty }\frac{\limfunc{essinf}\limits_{t<\tau <\infty
}\varphi _{1}(x_{0},\tau )\tau ^{\frac{\gamma }{p}}}{t^{\frac{\gamma }{p}+1}}%
dt\leq C\, \varphi _{2}(x_{0},r),  \label{e37}
\end{equation}%
and for $1<p<s$ the pair $(\varphi _{1},\varphi _{2})$ satisfies the
condition%
\begin{equation}
\int \limits_{r}^{\infty }\frac{\limfunc{essinf}\limits_{t<\tau <\infty
}\varphi _{1}(x_{0},\tau )\tau ^{\frac{\gamma }{p}}}{t^{\frac{\gamma }{p}-%
\frac{\gamma }{s}+1}}dt\leq C\, \varphi _{2}(x_{0},r)r^{\frac{\gamma }{s}},
\label{317}
\end{equation}%
where $C$ does not depend on $r$.

Then the operator $T_{\Omega }^{P}$ is bounded from $LM_{p,\varphi
_{1},P}^{\{x_{0}\}}$ to $LM_{p,\varphi _{2},P}^{\{x_{0}\}}$ for $p>1$ and
from $LM_{1,\varphi _{1},P}^{\{x_{0}\}}$ to $WLM_{1,\varphi
_{2},P}^{\{x_{0}\}}$ for $p=1$. Moreover, we have for $p>1$%
\begin{equation}
\left \Vert T_{\Omega }^{P}f\right \Vert _{LM_{p,\varphi
_{2},P}^{\{x_{0}\}}}\lesssim \left \Vert f\right \Vert _{LM_{p,\varphi
_{1},P}^{\{x_{0}\}}},  \label{3-1}
\end{equation}%
and for $p=1$%
\begin{equation}
\left \Vert T_{\Omega }^{P}f\right \Vert _{WLM_{1,\varphi
_{2},P}^{\{x_{0}\}}}\lesssim \left \Vert f\right \Vert _{LM_{1,\varphi
_{1},P}^{\{x_{0}\}}}.  \label{3-2}
\end{equation}
\end{theorem}

\begin{proof}
Let $1<p<\infty $ and $s^{\prime }\leq p$. By Theorem \ref{teo2*} and
Theorem \ref{teo5} with $v_{2}\left( r\right) =\varphi _{2}\left(
x_{0},r\right) ^{-1}$, $v_{1}=\varphi _{1}\left( x_{0},r\right) ^{-1}r^{-%
\frac{\gamma }{p}}$, $w\left( r\right) =r^{-\frac{\gamma }{p}-1}$ and $%
g\left( r\right) =\left \Vert f\right \Vert _{L_{p}\left( E\left(
x_{0},r\right) \right) }$, we have%
\begin{eqnarray*}
\left \Vert T_{\Omega }^{P}f\right \Vert _{LM_{p,\varphi
_{2},P}^{\{x_{0}\}}} &\lesssim &\sup_{r>0}\varphi _{2}\left( x_{0},r\right)
^{-1}\int \limits_{r}^{\infty }\left \Vert f\right \Vert _{L_{p}\left(
E\left( x_{0},t\right) \right) }\frac{dt}{t^{\frac{\gamma }{p}+1}} \\
&\lesssim &\sup_{r>0}\varphi _{1}\left( x_{0},r\right) ^{-1}r^{-\frac{\gamma 
}{p}}\left \Vert f\right \Vert _{L_{p}\left( E\left( x_{0},r\right) \right)
}=\left \Vert f\right \Vert _{_{LM_{p,\varphi _{1},P}^{\{x_{0}\}}}},
\end{eqnarray*}%
where the condition (\ref{3*}) is equivalent to (\ref{e37}), then we obtain (%
\ref{3-1}).

Let $1<p<s$. By Theorem \ref{teo2*} and Theorem \ref{teo5} with $v_{2}\left(
r\right) =\varphi _{2}\left( x_{0},r\right) ^{-1}$, $v_{1}=\varphi
_{1}\left( x_{0},r\right) ^{-1}r^{-\frac{\gamma }{p}+\frac{\gamma }{s}}$, $%
w\left( r\right) =r^{-\frac{\gamma }{p}+\frac{\gamma }{s}-1}$ and $g\left(
r\right) =\left \Vert f\right \Vert _{L_{p}\left( E\left( x_{0},r\right)
\right) }$, we have%
\begin{eqnarray*}
\left \Vert T_{\Omega }^{P}f\right \Vert _{_{LM_{p,\varphi
_{2},P}^{\{x_{0}\}}}} &\lesssim &\sup_{r>0}\varphi _{2}\left( x_{0},r\right)
^{-1}r^{-\frac{\gamma }{s}}\int \limits_{r}^{\infty }\left \Vert f\right
\Vert _{L_{p}\left( E\left( x_{0},t\right) \right) }\frac{dt}{t^{\frac{%
\gamma }{p}-\frac{\gamma }{s}+1}} \\
&\lesssim &\sup_{r>0}\varphi _{1}\left( x_{0},r\right) ^{-1}r^{-\frac{\gamma 
}{p}}\left \Vert f\right \Vert _{L_{p}\left( E\left( x_{0},r\right) \right)
}=\left \Vert f\right \Vert _{LM_{p,\varphi _{1},P}^{\{x_{0}\}}},
\end{eqnarray*}%
where the condition (\ref{3*}) is equivalent to (\ref{317}). Thus, we obtain
(\ref{3-1}).

Also, for $p=1$ we have%
\begin{eqnarray*}
\left \Vert T_{\Omega }^{P}f\right \Vert _{WLM_{1,\varphi
_{2},P}^{\{x_{0}\}}} &\lesssim &\sup_{r>0}\varphi _{2}\left( x_{0},r\right)
^{-1}\int \limits_{r}^{\infty }\left \Vert f\right \Vert _{L_{1}\left(
E\left( x_{0},t\right) \right) }\frac{dt}{t^{\gamma +1}} \\
&\lesssim &\sup_{r>0}\varphi _{1}\left( x_{0},r\right) ^{-1}r^{-\gamma
}\left \Vert f\right \Vert _{L_{1}\left( E\left( x_{0},r\right) \right)
}=\left \Vert f\right \Vert _{LM_{1,\varphi _{1},P}^{\{x_{0}\}}}.
\end{eqnarray*}%
Hence, the proof is completed.
\end{proof}

In the case of $s=\infty $ from Theorem \ref{teo9}, we get

\begin{corollary}
Let $x_{0}\in {\mathbb{R}^{n}}$, $1\leq p<\infty $ and the pair $(\varphi
_{1},\varphi _{2})$ satisfies condition (\ref{e37}). Then the operators $%
M^{P}$ and $\overline{T}^{P}$ are bounded from $LM_{p,\varphi
_{1},P}^{\{x_{0}\}}$ to $LM_{p,\varphi _{2},P}^{\{x_{0}\}}$ for $p>1$ and
from $LM_{1,\varphi _{1},P}^{\{x_{0}\}}$ to $WLM_{1,\varphi
_{2},P}^{\{x_{0}\}}$ for $p=1$.
\end{corollary}

\begin{corollary}
Let $x_{0}\in {\mathbb{R}^{n}}$, $1\leq p<\infty $ and $\Omega \in
L_{s}(S^{n-1})$, $1<s\leq \infty $, be $A_{t}$-homogeneous of degree zero.
For $s^{\prime }\leq p$, $p\neq 1$, the pair $\left( \varphi _{1},\varphi
_{2}\right) $ satisfies condition (\ref{e37}) and for $1<$ $p<s$ the pair $%
\left( \varphi _{1},\varphi _{2}\right) $ satisfies condition (\ref{317}).
Then the operators $M_{\Omega }^{P}$ and $\overline{T}_{\Omega }^{P}$ are
bounded from $LM_{p,\varphi _{1},P}^{\{x_{0}\}}$ to $LM_{p,\varphi
_{2},P}^{\{x_{0}\}}$ for $p>1$ and from $LM_{1,\varphi _{1},P}^{\{x_{0}\}}$
to $WLM_{1,\varphi _{2},P}^{\{x_{0}\}}$ for $p=1$.
\end{corollary}

\begin{corollary}
\label{Corollary 1}Let $x_{0}\in {\mathbb{R}^{n}}$, $1\leq p<\infty $ and $%
\Omega \in L_{s}(S^{n-1})$, $1<s\leq \infty $, is homogeneous of degree
zero. Let $T_{\Omega }^{P}$ be a parabolic sublinear operator satisfying
condition (\ref{e1}), bounded on $L_{p}({\mathbb{R}^{n}})$ for $p>1$, and
bounded from $L_{1}({\mathbb{R}^{n}})$ to $WL_{1}({\mathbb{R}^{n}})$. Let
also, for $s^{\prime }\leq p$, $p\neq 1$, the pair $\left( \varphi
_{1},\varphi _{2}\right) $ satisfies the condition 
\begin{equation*}
\int \limits_{r}^{\infty }\frac{\limfunc{essinf}\limits_{t<\tau <\infty
}\varphi _{1}(x_{0},\tau )\tau ^{\frac{n}{p}}}{t^{\frac{n}{p}+1}}dt\leq C\,
\varphi _{2}(x_{0},r),
\end{equation*}%
and for $1<p<s$ the pair $\left( \varphi _{1},\varphi _{2}\right) $
satisfies the condition%
\begin{equation*}
\int \limits_{r}^{\infty }\frac{\limfunc{essinf}\limits_{t<\tau <\infty
}\varphi _{1}(x_{0},\tau )\tau ^{\frac{n}{p}}}{t^{\frac{n}{p}-\frac{n}{s}+1}}%
dt\leq C\, \varphi _{2}(x_{0},r)r^{\frac{n}{s}},
\end{equation*}%
where $C$ does not depend on $r$.

Then the operator $T_{\Omega }^{P}$ is bounded from $LM_{p,\varphi
_{1},P}^{\{x_{0}\}}$ to $LM_{p,\varphi _{2},P}^{\{x_{0}\}}$ for $p>1$ and
from $LM_{1,\varphi _{1},P}^{\{x_{0}\}}$ to $WLM_{1,\varphi
_{2},P}^{\{x_{0}\}}$ for $p=1$. Moreover, we have for $p>1$%
\begin{equation*}
\left \Vert T_{\Omega }^{P}f\right \Vert _{LM_{p,\varphi
_{2},P}^{\{x_{0}\}}}\lesssim \left \Vert f\right \Vert _{LM_{p,\varphi
_{1},P}^{\{x_{0}\}}},
\end{equation*}%
and for $p=1$%
\begin{equation*}
\left \Vert T_{\Omega }^{P}f\right \Vert _{WLM_{1,\varphi
_{2},P}^{\{x_{0}\}}}\lesssim \left \Vert f\right \Vert _{LM_{1,\varphi
_{1},P}^{\{x_{0}\}}}.
\end{equation*}
\end{corollary}

\begin{remark}
Note that, in the case of $P=I$ Corollary \ref{Corollary 1} has been proved
in \cite{BGGS, Gurbuz1}. Also, in the case of $P=I$ and $s=\infty $
Corollary \ref{Corollary 1} has been proved in \cite{BGGS, Gurbuz1}.
\end{remark}

\begin{corollary}
\label{Corollary 2}Let $1\leq p<\infty $ and $\Omega \in L_{s}(S^{n-1})$, $%
1<s\leq \infty $, be $A_{t}$-homogeneous of degree zero. Let $T_{\Omega
}^{P} $ be a parabolic sublinear operator satisfying condition (\ref{e1}),
bounded on $L_{p}({\mathbb{R}^{n}})$ for $p>1$, and bounded from $L_{1}({%
\mathbb{R}^{n}})$ to $WL_{1}({\mathbb{R}^{n}})$. Let also, for $s^{\prime
}\leq p$, $p\neq 1$, the pair $(\varphi _{1},\varphi _{2})$ satisfies the
condition%
\begin{equation}
\int \limits_{r}^{\infty }\frac{\limfunc{essinf}\limits_{t<\tau <\infty
}\varphi _{1}(x,\tau )\tau ^{\frac{\gamma }{p}}}{t^{\frac{\gamma }{p}+1}}%
dt\leq C\, \varphi _{2}(x,r),  \label{e37*}
\end{equation}%
and for $1<p<s$ the pair $(\varphi _{1},\varphi _{2})$ satisfies the
condition%
\begin{equation}
\int \limits_{r}^{\infty }\frac{\limfunc{essinf}\limits_{t<\tau <\infty
}\varphi _{1}(x,\tau )\tau ^{\frac{\gamma }{p}}}{t^{\frac{\gamma }{p}-\frac{%
\gamma }{s}+1}}dt\leq C\, \varphi _{2}(x,r)r^{\frac{\gamma }{s}},
\label{e38}
\end{equation}%
where $C$ does not depend on $x$ and $r$.

Then the operator $T_{\Omega }^{P}$ is bounded from $M_{p,\varphi _{1},P}$
to $M_{p,\varphi _{2},P}$ for $p>1$ and from $M_{1,\varphi _{1},P}$ to $%
WM_{1,\varphi _{2},P}$ for $p=1$. Moreover, we have for $p>1$%
\begin{equation*}
\left \Vert T_{\Omega }^{P}f\right \Vert _{M_{p,\varphi _{2},P}}\lesssim
\left \Vert f\right \Vert _{M_{p,\varphi _{1},P}},
\end{equation*}%
and for $p=1$%
\begin{equation*}
\left \Vert T_{\Omega }^{P}f\right \Vert _{WM_{1,\varphi _{2},P}}\lesssim
\left \Vert f\right \Vert _{M_{1,\varphi _{1},P}}.
\end{equation*}
\end{corollary}

In the case of $s=\infty $ from Corollary \ref{Corollary 2}, we get

\begin{corollary}
Let $1\leq p<\infty $ and the pair $(\varphi _{1},\varphi _{2})$ satisfies
condition (\ref{e37*}). Then the operators $M^{P}$ and $\overline{T}^{P}$
are bounded from $M_{p,\varphi _{1},P}$ to $M_{p,\varphi _{2},P}$ for $p>1$
and from $M_{1,\varphi _{1},P}$ to $WM_{1,\varphi _{2},P}$ for $p=1$.
\end{corollary}

\begin{corollary}
Let $1\leq p<\infty $ and $\Omega \in L_{s}\left( S^{n-1}\right) $, $1<s\leq
\infty $, be $A_{t}$-homogeneous of degree zero. Let also, for $s^{\prime
}\leq p$, $p\neq 1$, the pair $\left( \varphi _{1},\varphi _{2}\right) $
satisfies condition (\ref{e37*}) and for $1<$ $p<q$ the pair $\left( \varphi
_{1},\varphi _{2}\right) $ satisfies condition (\ref{e38}). Then the
operators $M_{\Omega }^{P}$ and $\overline{T}_{\Omega }^{P}$ are bounded
from $M_{p,\varphi _{1}}$ to $M_{p,\varphi _{2}}$ for $p>1$ and from $%
M_{1,\varphi _{1}}$to $WM_{1,\varphi _{2}}$ for $p=1$.
\end{corollary}

\begin{remark}
Condition (\ref{e37*}) in Corollary \ref{Corollary 2} is weaker than
condition (\ref{e35}) in Theorem \ref{teo4}. Indeed, if condition (\ref{e35}%
) holds, then%
\begin{equation*}
\int \limits_{r}^{\infty }\frac{\limfunc{essinf}\limits_{t<\tau <\infty
}\varphi _{1}(x,\tau )\tau ^{\frac{\gamma }{p}}}{t^{\frac{\gamma }{p}+1}}%
dt\leq \int \limits_{r}^{\infty }\varphi _{1}(x,t)\frac{dt}{t},\qquad r\in
\left( 0,\infty \right) ,
\end{equation*}%
so condition (\ref{e37*}) holds.

On the other hand, the functions%
\begin{equation*}
\varphi _{1}(r)=\frac{1}{\chi _{\left( 1,\infty \right) }\left( r\right) r^{%
\frac{\gamma }{p}-\beta }},\qquad \varphi _{2}(r)=r^{-\frac{\gamma }{p}%
}\left( 1+r^{\beta }\right) ,\qquad 0<\beta <\frac{\gamma }{p}
\end{equation*}%
satisfy condition (\ref{e37*}) but do not satisfy condition (\ref{e35}) (see 
\cite{GULAKShIEOT2012, Karaman}).
\end{remark}

\begin{corollary}
\label{Corollary 3}Let $1\leq p<\infty $ and $\Omega \in L_{s}(S^{n-1})$, $%
1<s\leq \infty $, be homogeneous of degree zero. Let $T_{\Omega }^{P}$ be a
parabolic sublinear operator satisfying condition (\ref{e1}), bounded on $%
L_{p}({\mathbb{R}^{n}})$ for $p>1$, and bounded from $L_{1}({\mathbb{R}^{n}}%
) $ to $WL_{1}({\mathbb{R}^{n}})$. Let also, for $s^{\prime }\leq p$, $p\neq
1$, the pair $(\varphi _{1},\varphi _{2})$ satisfies the condition%
\begin{equation*}
\int \limits_{r}^{\infty }\frac{\limfunc{essinf}\limits_{t<\tau <\infty
}\varphi _{1}(x,\tau )\tau ^{\frac{n}{p}}}{t^{\frac{n}{p}+1}}dt\leq C\,
\varphi _{2}(x,r),
\end{equation*}%
and for $1<p<s$ the pair $(\varphi _{1},\varphi _{2})$ satisfies the
condition%
\begin{equation*}
\int \limits_{r}^{\infty }\frac{\limfunc{essinf}\limits_{t<\tau <\infty
}\varphi _{1}(x,\tau )\tau ^{\frac{n}{p}}}{t^{\frac{n}{p}-\frac{n}{s}+1}}%
dt\leq C\, \varphi _{2}(x,r)r^{\frac{n}{s}},
\end{equation*}%
where $C$ does not depend on $x$ and $r$.

Then the operator $T_{\Omega }^{P}$ is bounded from $M_{p,\varphi _{1},P}$
to $M_{p,\varphi _{2},P}$ for $p>1$ and from $M_{1,\varphi _{1},P}$ to $%
WM_{1,\varphi _{2},P}$ for $p=1$. Moreover, we have for $p>1$%
\begin{equation*}
\left \Vert T_{\Omega }^{P}f\right \Vert _{M_{p,\varphi _{2},P}}\lesssim
\left \Vert f\right \Vert _{M_{p,\varphi _{1},P}},
\end{equation*}%
and for $p=1$%
\begin{equation*}
\left \Vert T_{\Omega }^{P}f\right \Vert _{WM_{1,\varphi _{2},P}}\lesssim
\left \Vert f\right \Vert _{M_{1,\varphi _{1},P}}.
\end{equation*}
\end{corollary}

\begin{remark}
Note that, in the case of $P=I$ Corollary \ref{Corollary 3} has been proved
in \cite{BGGS, Gurbuz1, Gurbuz2}. Also, in the case of $P=I$ and $s=\infty $
Corollary \ref{Corollary 3} has been proved in \cite{BGGS, GULAKShIEOT2012,
Gurbuz1, Gurbuz2, Karaman}.
\end{remark}

\section{commutators of parabolic linear operators with rough kernel
generated by parabolic Calder\'{o}n-Zygmund operators and parabolic local
Campanato functions on the spaces $LM_{p,\protect \varphi ,P}^{\{x_{0}\}}$}

In this section, we will prove the boundedness of the operators $%
[b,T_{\Omega }^{P}]$ with $b\in LC_{p_{2},\lambda ,P}^{\left \{
x_{0}\right
\} }$ on the parabolic generalized local Morrey spaces $%
LM_{p,\varphi ,P}^{\{x_{0}\}}$ by using the following weighted Hardy operator

\begin{equation*}
H_{\omega }g(r):=\int \limits_{r}^{\infty }\left( 1+\ln \frac{t}{r}\right)
g(t)\omega (t)dt,\qquad r\in \left( 0,\infty \right) ,
\end{equation*}%
where $\omega $ is a weight function.

Let $T$ be a linear operator. For a locally integrable function $b$ on ${%
\mathbb{R}^{n}}$, we define the commutator $[b,T]$ by 
\begin{equation*}
\lbrack b,T]f(x)=b(x)\,Tf(x)-T(bf)(x)
\end{equation*}%
for any suitable function $f$. Let $\overline{T}$ be a Calder\'{o}n--Zygmund
operator. A well known result of Coifman et al. \cite{CRW} states that when $%
K\left( x\right) =\frac{\Omega \left( x^{\prime }\right) }{\left \vert
x\right \vert ^{n}}$ and $\Omega $ is smooth, the commutator $[b,\overline{T}%
]f=b\, \overline{T}f-\overline{T}(bf)$ is bounded on $L_{p}({\mathbb{R}^{n}}%
) $, $1<p<\infty $, if and only if $b\in BMO({\mathbb{R}^{n}})$.

Since $BMO({\mathbb{R}^{n}})\subset \bigcap \limits_{p>1}LC_{p,P}^{\left \{
x_{0}\right \} }({\mathbb{R}^{n}})$, if we only assume $b\in
LC_{p,P}^{\left
\{ x_{0}\right \} }({\mathbb{R}^{n}})$, or more generally $%
b\in LC_{p,\lambda ,P}^{\left \{ x_{0}\right \} }({\mathbb{R}^{n}})$, then $%
[b,\overline{T}]$ may not be a bounded operator on $L_{p}({\mathbb{R}^{n}})$%
, $1<p<\infty $. However, it has some boundedness properties on other
spaces. As a matter of fact, Grafakos et al. \cite{GraLiYang} have
considered the commutator with $b\in LC_{p,I}^{\left \{ x_{0}\right \} }({%
\mathbb{R}^{n}})$ on Herz spaces for the first time. Morever, in \cite{BGGS}%
, \cite{Lu1}, \cite{Gurbuz1} and \cite{TaoShi}, they have considered the
commutators with $b\in LC_{p,\lambda ,I}^{\left \{ x_{0}\right \} }({\mathbb{%
R}^{n}})$. The commutator of Calder\'{o}n--Zygmund operators plays an
important role in studying the regularity of solutions of elliptic partial
differential equations of second order (see, for example, \cite{ChFraL1,
ChFraL2, FazRag2}). The boundedness of the commutator has been generalized
to other contexts and important applications to some non-linear PDEs have
been given by Coifman et al. \cite{CLMS}.

We introduce the parabolic local Campanato space $LC_{p,\lambda
,P}^{\left
\{ x_{0}\right \} }$ following the known ideas of defining local
Campanato space (see \cite{BGGS, Gurbuz1} etc).

\begin{definition}
Let $1\leq p<\infty $ and $0\leq \lambda <\frac{1}{\gamma }$. A function $%
f\in L_{p}^{loc}\left( {\mathbb{R}^{n}}\right) $ is said to belong to the $%
LC_{p,\lambda ,P}^{\left \{ x_{0}\right \} }\left( {\mathbb{R}^{n}}\right) $
(parabolic local Campanato space), if%
\begin{equation}
\left \Vert f\right \Vert _{LC_{p,\lambda ,P}^{\left \{ x_{0}\right \}
}}=\sup_{r>0}\left( \frac{1}{\left \vert E\left( x_{0},r\right) \right \vert
^{1+\lambda p}}\int \limits_{E\left( x_{0},r\right) }\left \vert f\left(
y\right) -f_{E\left( x_{0},r\right) }\right \vert ^{p}dy\right) ^{\frac{1}{p}%
}<\infty ,  \label{e51}
\end{equation}%
where%
\begin{equation*}
f_{E\left( x_{0},r\right) }=\frac{1}{\left \vert E\left( x_{0},r\right)
\right \vert }\int \limits_{E\left( x_{0},r\right) }f\left( y\right) dy.
\end{equation*}

Define%
\begin{equation*}
LC_{p,\lambda ,P}^{\left \{ x_{0}\right \} }\left( {\mathbb{R}^{n}}\right)
=\left \{ f\in L_{p}^{loc}\left( {\mathbb{R}^{n}}\right) :\left \Vert
f\right \Vert _{LC_{p,\lambda ,P}^{\left \{ x_{0}\right \} }}<\infty \right
\} .
\end{equation*}
\end{definition}

\begin{remark}
If two functions which differ by a constant are regarded as a function in
the space $LC_{p,\lambda ,P}^{\left \{ x_{0}\right \} }\left( {\mathbb{R}^{n}%
}\right) $, then $LC_{p,\lambda ,P}^{\left \{ x_{0}\right \} }\left( {%
\mathbb{R}^{n}}\right) $ becomes a Banach space. The space $LC_{p,\lambda
,P}^{\left \{ x_{0}\right \} }\left( {\mathbb{R}^{n}}\right) $ when $\lambda
=0$ is just the $LC_{p,P}^{\left \{ x_{0}\right \} }({\mathbb{R}^{n}})$.
Apparently, (\ref{e51}) is equivalent to the following condition:%
\begin{equation*}
\sup_{r>0}\inf_{c\in 
\mathbb{C}
}\left( \frac{1}{\left \vert E\left( x_{0},r\right) \right \vert ^{1+\lambda
p}}\int \limits_{E\left( x_{0},r\right) }\left \vert f\left( y\right)
-c\right \vert ^{p}dy\right) ^{\frac{1}{p}}<\infty .
\end{equation*}
\end{remark}

In \cite{LuYang1}, Lu and Yang has introduced the central BMO space $%
CBMO_{p}({\mathbb{R}^{n}})=LC_{p,0,I}^{\{0\}}({\mathbb{R}^{n}})$. Also the
space $CBMO^{\{x_{0}\}}({\mathbb{R}^{n}})=LC_{1,0,I}^{\{x_{0}\}}({\mathbb{R}%
^{n}})$ has been considered in other denotes in \cite{Rzaev}. The space $%
LC_{p,P}^{\left \{ x_{0}\right \} }({\mathbb{R}^{n}})$ can be regarded as a
local version of $BMO({\mathbb{R}^{n}})$, the space of parabolic bounded
mean oscillation, at the origin. But, they have quite different properties.
The classical John-Nirenberg inequality shows that functions in $BMO({%
\mathbb{R}^{n}})$ are locally exponentially integrable. This implies that,
for any $1\leq p<\infty $, the functions in $BMO({\mathbb{R}^{n}})$
(parabolic $BMO$) can be described by means of the condition:%
\begin{equation*}
\sup_{x\in {\mathbb{R}^{n},r>0}}\left( \frac{1}{|E\left( x,r\right) |}\dint
\limits_{E\left( x,r\right) }|f(y)-f_{E\left( x,r\right) }|^{p}dy\right) ^{%
\frac{1}{p}}<\infty ,
\end{equation*}%
where $B$ denotes an arbitrary ball in ${\mathbb{R}^{n}}$. However, the
space $LC_{p,P}^{\left \{ x_{0}\right \} }({\mathbb{R}^{n}})$ depends on $p$%
. If $p_{1}<p_{2}$, then $LC_{p_{2},P}^{\left \{ x_{0}\right \} }({\mathbb{R}%
^{n}})\subsetneqq LC_{p_{1},P}^{\left \{ x_{0}\right \} }({\mathbb{R}^{n}})$%
. Therefore, there is no analogy of the famous John-Nirenberg inequality of $%
BMO({\mathbb{R}^{n}})$ for the space $LC_{p,P}^{\left \{ x_{0}\right \} }({%
\mathbb{R}^{n}})$. One can imagine that the behavior of $LC_{p,P}^{\left \{
x_{0}\right \} }({\mathbb{R}^{n}})$ may be quite different from that of $BMO(%
{\mathbb{R}^{n}})$.

\begin{theorem}
\label{teo5*}\cite{BGGS, Gurbuz1} Let $v_{1}$, $v_{2}$ and $\omega $ be
weigths on $(0,\infty )$ and $v_{1}\left( t\right) $ be bounded outside a
neighbourhood of the origin. The inequality 
\begin{equation}
\limfunc{esssup}\limits_{r>0}v_{2}(r)H_{\omega }g(r)\leq C\limfunc{esssup}%
\limits_{r>0}v_{1}(r)g(r)  \label{e52}
\end{equation}%
holds for some $C>0$ for all non-negative and non-decreasing functions $g$
on $(0,\infty )$ if and only if 
\begin{equation}
B:=\sup \limits_{r>0}v_{2}(r)\int \limits_{r}^{\infty }\left( 1+\ln \frac{t}{%
r}\right) \frac{\omega (t)dt}{\limfunc{esssup}\limits_{t<s<\infty }v_{1}(s)}%
<\infty .  \label{3}
\end{equation}%
Moreover, the value $C=B$ is the best constant for (\ref{e52}).
\end{theorem}

\begin{remark}
In (\ref{e52}) and (\ref{3}) it is assumed that $\frac{1}{\infty }=0$ and $%
0.\infty =0$.
\end{remark}

\begin{lemma}
\label{Lemma 4}Let $b$ be function in $LC_{p,\lambda ,P}^{\left \{
x_{0}\right \} }\left( 
\mathbb{R}
^{n}\right) $, $1\leq p<\infty $, $0\leq \lambda <\frac{1}{\gamma }$ and $%
r_{1}$, $r_{2}>0$. Then%
\begin{equation}
\left( \frac{1}{\left \vert E\left( x_{0},r_{1}\right) \right \vert
^{1+\lambda p}}\dint \limits_{E\left( x_{0},r_{1}\right) }\left \vert
b\left( y\right) -b_{E\left( x_{0},r_{2}\right) }\right \vert ^{p}dy\right)
^{\frac{1}{p}}\leq C\left( 1+\ln \frac{r_{1}}{r_{2}}\right) \left \Vert
b\right \Vert _{LC_{p,\lambda ,P}^{\left \{ x_{0}\right \} }},  \label{a}
\end{equation}%
where $C>0$ is independent of $b$, $r_{1}$ and $r_{2}$.

From this inequality $\left( \text{\ref{a}}\right) $, we have%
\begin{equation}
\left \vert b_{E\left( x_{0},r_{1}\right) }-b_{E\left( x_{0},r_{2}\right)
}\right \vert \leq C\left( 1+\ln \frac{r_{1}}{r_{2}}\right) \left \vert
E\left( x_{0},r_{1}\right) \right \vert ^{\lambda }\left \Vert b\right \Vert
_{LC_{p,\lambda ,P}^{\left \{ x_{0}\right \} }},  \label{b}
\end{equation}

and it is easy to see that%
\begin{equation}
\left \Vert b-b_{E}\right \Vert _{L_{p}\left( E\right) }\leq C\left( 1+\ln 
\frac{r_{1}}{r_{2}}\right) r^{\frac{\gamma }{p}+\gamma \lambda }\left \Vert
b\right \Vert _{LC_{p,\lambda ,P}^{\left \{ x_{0}\right \} }}.  \label{c}
\end{equation}
\end{lemma}

In \cite{DingYZ} the following statements have been proved for the parabolic
commutators of parabolic singular integral operators with rough kernel $%
\overline{T}_{\Omega }^{P}$, containing the result in \cite{Miz, Nakai1,
Nakai2}.

\begin{theorem}
\textit{Suppose that }$\Omega \in L_{s}(S^{n-1})$\textit{, }$1<s\leq \infty $%
\textit{, is }$A_{t}$-\textit{homogeneous of degree zero }and $b\in BMO({%
\mathbb{R}^{n}})$\textit{. Let} $1\leq s^{\prime }<p<\infty \left( s^{\prime
}=\frac{s}{s-1}\right) $ and\textit{\ }$\varphi (x,r)$ satisfies the
conditions (\ref{e32}) and (\ref{e33x}). \textit{If the commutator operator }%
$[b,\overline{T}_{\Omega }^{P}]$\textit{\ is bounded on }$L_{p}\left( {%
\mathbb{R}^{n}}\right) $, then the operator $[b,\overline{T}_{\Omega }^{P}]$
is bounded on $M_{p,\varphi ,P}$.
\end{theorem}

\begin{theorem}
Let $1<p<\infty $, $b\in BMO({\mathbb{R}^{n}})$ and $\varphi (x,t)$
satisfies conditions (\ref{e32}) and (\ref{e33x}). Then the operators $%
M_{b}^{P}$ and $[b,\overline{T}^{P}]$ are bounded on $M_{p,\varphi ,P}$.
\end{theorem}

As in the proof of Theorem \ref{teo9}, it suffices to prove the following
Theorem \ref{teo2} (our main result).

\begin{theorem}
\label{teo2}Let $x_{0}\in {\mathbb{R}^{n}}$, $1<p<\infty $ and $\Omega \in
L_{s}(S^{n-1})$, $1<s\leq \infty $, be $A_{t}$-homogeneous of degree zero.
Let $T_{\Omega }^{P}$ be a parabolic linear operator satisfying condition (%
\ref{e1}), bounded on $L_{p}({\mathbb{R}^{n}})$ for $1<p<\infty $. Let also, 
$b\in LC_{p_{2},\lambda ,P}^{\left \{ x_{0}\right \} }\left( 
\mathbb{R}
^{n}\right) $, $0\leq \lambda <\frac{1}{\gamma }$ and $\frac{1}{p}=\frac{1}{%
p_{1}}+\frac{1}{p_{2}}$.

\textit{Then, for }$s^{\prime }\leq p$\textit{, the inequality }%
\begin{equation*}
\left \Vert \lbrack b,T_{\Omega }^{P}]f\right \Vert _{L_{p}\left( E\left(
x_{0},r\right) \right) }\lesssim \Vert b\Vert _{LC_{p_{2},\lambda
,P}^{\{x_{0}\}}}r^{\frac{\gamma }{p}}\int \limits_{2kr}^{\infty }\left(
1+\ln \frac{t}{r}\right) t^{\gamma \lambda -\frac{\gamma }{p_{1}}-1}\left
\Vert f\right \Vert _{L_{p_{1}}\left( E\left( x_{0},t\right) \right) }dt
\end{equation*}%
holds for any ellipsoid $E\left( x_{0},r\right) $ and for all $f\in
L_{p_{1}}^{loc}\left( {\mathbb{R}^{n}}\right) $.

Also, for $p_{1}<s$, the inequality$\,$%
\begin{equation*}
\left \Vert \lbrack b,T_{\Omega }^{P}]f\right \Vert _{L_{p}\left( E\left(
x_{0},r\right) \right) }\lesssim \Vert b\Vert _{LC_{p_{2},\lambda
,P}^{\{x_{0}\}}}r^{\frac{\gamma }{p}-\frac{\gamma }{s}}\int
\limits_{2kr}^{\infty }\left( 1+\ln \frac{t}{r}\right) t^{\gamma \lambda -%
\frac{\gamma }{p_{1}}+\frac{\gamma }{s}-1}\left \Vert f\right \Vert
_{L_{p_{1}}\left( E\left( x_{0},t\right) \right) }dt
\end{equation*}%
holds for any ellipsoid $E\left( x_{0},r\right) $ and for all $f\in
L_{p_{1}}^{loc}\left( {\mathbb{R}^{n}}\right) $.
\end{theorem}

\begin{proof}
Let $1<p<\infty $, $b\in LC_{p_{2},\lambda ,P}^{\left \{ x_{0}\right \} }({%
\mathbb{R}^{n}})$ and $\frac{1}{p}=\frac{1}{p_{1}}+\frac{1}{p_{2}}$. Set $%
E=E\left( x_{0},r\right) $ for the parabolic ball (ellipsoid) centered at $%
x_{0}$ and of radius $r$ and $2kE=E\left( x_{0},2kr\right) $. We represent $%
f $ as%
\begin{equation*}
f=f_{1}+f_{2},\qquad \text{\ }f_{1}\left( y\right) =f\left( y\right) \chi
_{2kE}\left( y\right) ,\qquad \text{\ }f_{2}\left( y\right) =f\left(
y\right) \chi _{\left( 2kE\right) ^{C}}\left( y\right) ,\qquad r>0
\end{equation*}%
and have%
\begin{align*}
\lbrack b,T_{\Omega }^{P}]f\left( x\right) & =\left( b\left( x\right)
-b_{E}\right) T_{\Omega }^{P}f_{1}\left( x\right) -T_{\Omega }^{P}\left(
\left( b\left( \cdot \right) -b_{E}\right) f_{1}\right) \left( x\right) \\
& +\left( b\left( x\right) -b_{E}\right) T_{\Omega }^{P}f_{2}\left( x\right)
-T_{\Omega }^{P}\left( \left( b\left( \cdot \right) -b_{E}\right)
f_{2}\right) \left( x\right) \\
& \equiv J_{1}+J_{2}+J_{3}+J_{4}.
\end{align*}

Hence we get%
\begin{equation*}
\left \Vert \lbrack b,T_{\Omega }^{P}]f\right \Vert _{L_{p}\left( E\right)
}\leq \left \Vert J_{1}\right \Vert _{L_{p}\left( E\right) }+\left \Vert
J_{2}\right \Vert _{L_{p}\left( E\right) }+\left \Vert J_{3}\right \Vert
_{L_{p}\left( E\right) }+\left \Vert J_{4}\right \Vert _{L_{p}\left(
E\right) }.
\end{equation*}%
By the H\"{o}lder's inequality, the boundedness of $T_{\Omega }^{P}$ on $%
L_{p_{1}}({\mathbb{R}^{n}})$ (see Theorem \ref{Teo-Guliyev}) it follows that:%
\begin{eqnarray*}
\left \Vert J_{1}\right \Vert _{L_{p}\left( E\right) } &\leq &\left \Vert
\left( b\left( \cdot \right) -b_{E}\right) I_{\Omega ,\alpha
}^{P}f_{1}\left( \cdot \right) \right \Vert _{L_{p}\left( {\mathbb{R}^{n}}%
\right) } \\
&\lesssim &\left \Vert \left( b\left( \cdot \right) -b_{E}\right) \right
\Vert _{L_{p_{2}}\left( {\mathbb{R}^{n}}\right) }\left \Vert I_{\Omega
,\alpha }^{P}f_{1}\left( \cdot \right) \right \Vert _{L_{p_{1}}\left( {%
\mathbb{R}^{n}}\right) } \\
&\lesssim &\left \Vert b\right \Vert _{LC_{p_{2},\lambda ,P}^{\left \{
x_{0}\right \} }}r^{\frac{\gamma }{p_{2}}+\gamma \lambda }\left \Vert
f_{1}\right \Vert _{L_{p_{1}}\left( {\mathbb{R}^{n}}\right) } \\
&=&\left \Vert b\right \Vert _{LC_{p_{2},\lambda ,P}^{\left \{ x_{0}\right
\} }}r^{\frac{\gamma }{p_{2}}+\frac{\gamma }{p_{1}}+\gamma \lambda }\left
\Vert f\right \Vert _{L_{p_{1}}\left( 2kE\right) }\int \limits_{2kr}^{\infty
}t^{-1-\frac{\gamma }{p_{1}}}dt \\
&\lesssim &\left \Vert b\right \Vert _{LC_{p_{2},\lambda ,P}^{\left \{
x_{0}\right \} }}r^{\frac{\gamma }{p}}\int \limits_{2kr}^{\infty }\left(
1+\ln \frac{t}{r}\right) t^{\gamma \lambda -\frac{\gamma }{p_{1}}-1}\left
\Vert f\right \Vert _{L_{p_{1}}\left( E\left( x_{0},t\right) \right) }dt.
\end{eqnarray*}

Using the the boundedness of $T_{\Omega }^{P}$ on $L_{p}({\mathbb{R}^{n}})$
(see Theorem \ref{Teo-Guliyev}), by the H\"{o}lder's inequality for $J_{2}$,
we have%
\begin{align*}
\left \Vert J_{2}\right \Vert _{L_{p}\left( E\right) }& \leq \left \Vert
I_{\Omega ,\alpha }^{P}\left( b\left( \cdot \right) -b_{E}\right)
f_{1}\right \Vert _{L_{p}\left( {\mathbb{R}^{n}}\right) } \\
& \lesssim \left \Vert \left( b\left( \cdot \right) -b_{E}\right)
f_{1}\right \Vert _{L_{p}\left( {\mathbb{R}^{n}}\right) } \\
& \lesssim \left \Vert b\left( \cdot \right) -b_{E}\right \Vert
_{L_{p_{2}}\left( {\mathbb{R}^{n}}\right) }\left \Vert f_{1}\right \Vert
_{L_{p_{1}}\left( {\mathbb{R}^{n}}\right) } \\
& \lesssim \left \Vert b\right \Vert _{LC_{p_{2},\lambda ,P}^{\left \{
x_{0}\right \} }}r^{\frac{\gamma }{p_{2}}+\frac{\gamma }{p_{1}}+\gamma
\lambda }\left \Vert f\right \Vert _{L_{p_{1}}\left( 2kE\right) }\int
\limits_{2kr}^{\infty }t^{-1-\frac{\gamma }{p_{1}}}dt \\
& \lesssim \left \Vert b\right \Vert _{LC_{p_{2},\lambda ,P}^{\left \{
x_{0}\right \} }}r^{\frac{\gamma }{p}}\int \limits_{2kr}^{\infty }\left(
1+\ln \frac{t}{r}\right) t^{\gamma \lambda -\frac{\gamma }{p_{1}}-1}\left
\Vert f\right \Vert _{L_{p_{1}}\left( E\left( x_{0},t\right) \right) }dt.
\end{align*}

For $J_{3}$, it is known that $x\in E$, $y\in \left( 2kE\right) ^{C}$, which
implies $\frac{1}{2k}\rho \left( x_{0}-y\right) \leq \rho \left( x-y\right)
\leq \frac{3k}{2}\rho \left( x_{0}-y\right) $.

When $s^{\prime }\leq p_{1}$, by the Fubini's theorem, the H\"{o}lder's
inequality and (\ref{e311}), we have%
\begin{align*}
\left \vert T_{\Omega }^{P}f_{2}\left( x\right) \right \vert & \leq
c_{0}\int \limits_{\left( 2kE\right) ^{C}}\left \vert \Omega \left(
x-y\right) \right \vert \frac{\left \vert f\left( y\right) \right \vert }{%
\rho \left( x_{0}-y\right) ^{\gamma }}dy \\
& \approx \int \limits_{2kr}^{\infty }\int \limits_{2kr<\rho \left(
x_{0}-y\right) <t}\left \vert \Omega \left( x-y\right) \right \vert \left
\vert f\left( y\right) \right \vert dyt^{-1-\gamma }dt \\
& \lesssim \int \limits_{2kr}^{\infty }\int \limits_{E\left( x_{0},t\right)
}\left \vert \Omega \left( x-y\right) \right \vert \left \vert f\left(
y\right) \right \vert dyt^{-1-\gamma }dt \\
& \lesssim \int \limits_{2kr}^{\infty }\left \Vert f\right \Vert
_{L_{p_{1}}\left( E\left( x_{0},t\right) \right) }\left \Vert \Omega \left(
x-\cdot \right) \right \Vert _{L_{s}\left( E\left( x_{0},t\right) \right)
}\left \vert E\left( x_{0},t\right) \right \vert ^{1-\frac{1}{p_{1}}-\frac{1%
}{s}}t^{-1-\gamma }dt \\
& \lesssim \int \limits_{2kr}^{\infty }\left \Vert f\right \Vert
_{L_{p_{1}}\left( E\left( x_{0},t\right) \right) }t^{-1-\frac{\gamma }{p_{1}}%
}dt.
\end{align*}

Hence, we get%
\begin{align*}
\left \Vert J_{3}\right \Vert _{L_{p}\left( E\right) }& \leq \left \Vert
\left( b\left( \cdot \right) -b_{B}\right) T_{\Omega }^{P}f_{2}\left( \cdot
\right) \right \Vert _{L_{p}\left( {\mathbb{R}^{n}}\right) } \\
& \lesssim \left \Vert \left( b\left( \cdot \right) -b_{E}\right) \right
\Vert _{L_{p}\left( {\mathbb{R}^{n}}\right) }\int \limits_{2kr}^{\infty
}t^{-1-\frac{\gamma }{p_{1}}}\left \Vert f\right \Vert _{L_{p_{1}}\left(
E\left( x_{0},t\right) \right) }dt \\
& \lesssim \left \Vert \left( b\left( \cdot \right) -b_{E}\right) \right
\Vert _{L_{p_{2}}\left( {\mathbb{R}^{n}}\right) }r^{\frac{\gamma }{p_{1}}%
}\int \limits_{2kr}^{\infty }t^{-1-\frac{\gamma }{p_{1}}}\left \Vert f\right
\Vert _{L_{p_{1}}\left( E\left( x_{0},t\right) \right) }dt \\
& \lesssim \left \Vert b\right \Vert _{LC_{p_{2},\lambda ,P}^{\left \{
x_{0}\right \} }}r^{\frac{\gamma }{p}+\gamma \lambda }\int
\limits_{2kr}^{\infty }\left( 1+\ln \frac{t}{r}\right) t^{-1-\frac{\gamma }{%
p_{1}}}\left \Vert f\right \Vert _{L_{p_{1}}\left( E\left( x_{0},t\right)
\right) }dt \\
& \lesssim \left \Vert b\right \Vert _{LC_{p_{2},\lambda ,P}^{\left \{
x_{0}\right \} }}r^{\frac{\gamma }{p}}\int \limits_{2kr}^{\infty }\left(
1+\ln \frac{t}{r}\right) t^{\gamma \lambda -\frac{\gamma }{p_{1}}-1}\left
\Vert f\right \Vert _{L_{p_{1}}\left( E\left( x_{0},t\right) \right) }dt.
\end{align*}

When $p_{1}<s$, by the Fubini's theorem, the Minkowski inequality, the H\"{o}%
lder's inequality and from (\ref{c}), (\ref{314}), we get%
\begin{align*}
\left \Vert J_{3}\right \Vert _{L_{p}\left( E\right) }& \leq \left( \int
\limits_{E}\left \vert \int \limits_{2kr}^{\infty }\int \limits_{E\left(
x_{0},t\right) }\left \vert f\left( y\right) \right \vert \left \vert
b\left( x\right) -b_{E}\right \vert \left \vert \Omega \left( x-y\right)
\right \vert dy\frac{dt}{t^{\gamma +1}}\right \vert ^{p}dx\right) ^{\frac{1}{%
p}} \\
& \leq \int \limits_{2kr}^{\infty }\int \limits_{E\left( x_{0},t\right)
}\left \vert f\left( y\right) \right \vert \left \Vert \left( b\left( \cdot
\right) -b_{E}\right) \Omega \left( \cdot -y\right) \right \Vert
_{L_{p}\left( E\right) }dy\frac{dt}{t^{\gamma +1}} \\
& \leq \int \limits_{2kr}^{\infty }\int \limits_{E\left( x_{0},t\right)
}\left \vert f\left( y\right) \right \vert \left \Vert b\left( \cdot \right)
-b_{E}\right \Vert _{L_{p_{2}}\left( E\right) }\left \Vert \Omega \left(
\cdot -y\right) \right \Vert _{L_{p_{1}}\left( E\right) }dy\frac{dt}{%
t^{\gamma +1}} \\
& \lesssim \left \Vert b\right \Vert _{LC_{p_{2},\lambda ,P}^{\left \{
x_{0}\right \} }}r^{\frac{\gamma }{p_{2}}+\gamma \lambda }\left \vert
E\right \vert ^{\frac{1}{p_{1}}-\frac{1}{s}}\int \limits_{2kr}^{\infty }\int
\limits_{E\left( x_{0},t\right) }\left \vert f\left( y\right) \right \vert
\left \Vert \Omega \left( \cdot -y\right) \right \Vert _{L_{s}\left(
E\right) }dy\frac{dt}{t^{\gamma +1}} \\
& \lesssim \left \Vert b\right \Vert _{LC_{p_{2},\lambda ,P}^{\left \{
x_{0}\right \} }}r^{\frac{\gamma }{p}-\frac{\gamma }{s}+\gamma \lambda }\int
\limits_{2kr}^{\infty }\left \Vert f\right \Vert _{L_{1}\left( E\left(
x_{0},t\right) \right) }\left \vert E\left( x_{0},\frac{3}{2}t\right) \right
\vert ^{\frac{1}{s}}\frac{dt}{t^{\gamma +1}} \\
& \lesssim \left \Vert b\right \Vert _{LC_{p_{2},\lambda ,P}^{\left \{
x_{0}\right \} }}r^{\frac{\gamma }{p}-\frac{\gamma }{s}+\gamma \lambda }\int
\limits_{2kr}^{\infty }\left( 1+\ln \frac{t}{r}\right) \left \Vert f\right
\Vert _{L_{p_{1}}\left( E\left( x_{0},t\right) \right) }\frac{dt}{t^{\frac{%
\gamma }{p_{1}}-\frac{\gamma }{s}+1}} \\
& \lesssim \left \Vert b\right \Vert _{LC_{p_{2},\lambda ,P}^{\left \{
x_{0}\right \} }}r^{\frac{\gamma }{p}-\frac{\gamma }{s}}\int
\limits_{2kr}^{\infty }\left( 1+\ln \frac{t}{r}\right) t^{\gamma \lambda -%
\frac{\gamma }{p_{1}}+\frac{\gamma }{s}-1}\left \Vert f\right \Vert
_{L_{p_{1}}\left( E\left( x_{0},t\right) \right) }dt.
\end{align*}

On the other hand, for $J_{4}$, when $s^{\prime }\leq p$, for $x\in E$, by
the Fubini's theorem, applying the H\"{o}lder's inequality and from (\ref{b}%
), (\ref{c}) (\ref{e311}) we have

$\left \vert I_{\Omega ,\alpha }^{P}\left( \left( b\left( \cdot \right)
-b_{B}\right) f_{2}\right) \left( x\right) \right \vert \lesssim \dint
\limits_{\left( 2kE\right) ^{C}}\left \vert b\left( y\right)
-b_{E}\right
\vert \left \vert \Omega \left( x-y\right) \right \vert \frac{%
\left \vert f\left( y\right) \right \vert }{\rho \left( x-y\right) ^{\gamma }%
}dy$

$\lesssim \dint \limits_{\left( 2kE\right) ^{C}}\left \vert b\left( y\right)
-b_{E}\right \vert \left \vert \Omega \left( x-y\right) \right \vert \frac{%
\left \vert f\left( y\right) \right \vert }{\rho \left( x_{0}-y\right)
^{\gamma }}dy$

$\approx \dint \limits_{2kr}^{\infty }\dint \limits_{2kr<\rho \left(
x_{0}-y\right) <t}\left \vert b\left( y\right) -b_{E}\right \vert
\left
\vert \Omega \left( x-y\right) \right \vert \left \vert f\left(
y\right) \right \vert dy\frac{dt}{t^{\gamma +1}}$

$\lesssim \dint \limits_{2kr}^{\infty }\dint \limits_{E\left( x_{0},t\right)
}\left \vert b\left( y\right) -b_{E\left( x_{0},t\right) }\right \vert
\left
\vert \Omega \left( x-y\right) \right \vert \left \vert f\left(
y\right) \right \vert dy\frac{dt}{t^{\gamma +1}}$

$+\dint \limits_{2kr}^{\infty }\left \vert b_{E\left( x_{0},r\right)
}-b_{E\left( x_{0},t\right) }\right \vert \dint \limits_{E\left(
x_{0},t\right) }\left \vert \Omega \left( x-y\right) \right \vert
\left
\vert f\left( y\right) \right \vert dy\frac{dt}{t^{\gamma +1}}$

$\lesssim \dint \limits_{2kr}^{\infty }\left \Vert \left( b\left( \cdot
\right) -b_{E\left( x_{0},t\right) }\right) f\right \Vert _{L_{p}\left(
E\left( x_{0},t\right) \right) }\left \Vert \Omega \left( \cdot -y\right)
\right \Vert _{L_{s}\left( E\left( x_{0},t\right) \right) }\left \vert
E\left( x_{0},t\right) \right \vert ^{1-\frac{1}{p}-\frac{1}{s}}\frac{dt}{%
t^{\gamma +1}}$

$+\dint \limits_{2kr}^{\infty }\left \vert b_{E\left( x_{0},r\right)
}-b_{E\left( x_{0},t\right) }\right \vert \left \Vert f\right \Vert
_{L_{p_{1}}\left( E\left( x_{0},t\right) \right) }\left \Vert \Omega \left(
\cdot -y\right) \right \Vert _{L_{s}\left( E\left( x_{0},t\right) \right)
}\left \vert E\left( x_{0},t\right) \right \vert ^{1-\frac{1}{p_{1}}-\frac{1%
}{s}}t^{-\gamma -1}dt$

$\lesssim \dint \limits_{2kr}^{\infty }\left \Vert \left( b\left( \cdot
\right) -b_{E\left( x_{0},t\right) }\right) \right \Vert _{L_{p_{2}}\left(
E\left( x_{0},t\right) \right) }\left \Vert f\right \Vert _{L_{p_{1}}\left(
E\left( x_{0},t\right) \right) }t^{-1-\frac{\gamma }{p_{1}}}dt$

$+\left \Vert b\right \Vert _{LC_{p_{2},\lambda ,P}^{\left \{ x_{0}\right \}
}}\dint \limits_{2kr}^{\infty }\left( 1+\ln \frac{t}{r}\right) \left \Vert
f\right \Vert _{L_{p_{1}}\left( E\left( x_{0},t\right) \right) }t^{-1-\frac{%
\gamma }{p_{1}}+\gamma \lambda }dt$

$\lesssim \left \Vert b\right \Vert _{LC_{p_{2},\lambda ,P}^{\left \{
x_{0}\right \} }}\dint \limits_{2kr}^{\infty }\left( 1+\ln \frac{t}{r}%
\right) \left \Vert f\right \Vert _{L_{p_{1}}\left( E\left( x_{0},t\right)
\right) }t^{-1-\frac{\gamma }{p_{1}}+\gamma \lambda }dt.$

Then, we have%
\begin{eqnarray*}
\left \Vert J_{4}\right \Vert _{L_{p}\left( E\right) } &=&\left \Vert
I_{\Omega ,\alpha }^{P}\left( b\left( \cdot \right) -b_{E}\right)
f_{2}\right \Vert _{L_{p}\left( E\right) } \\
&\lesssim &\left \Vert b\right \Vert _{LC_{p_{2},\lambda ,P}^{\left \{
x_{0}\right \} }}r^{\frac{n}{p}}\dint \limits_{2kr}^{\infty }\left( 1+\ln 
\frac{t}{r}\right) t^{\gamma \lambda -\frac{\gamma }{p_{1}}-1}\left \Vert
f\right \Vert _{L_{p_{1}}\left( E\left( x_{0},t\right) \right) }dt.
\end{eqnarray*}

When $p_{1}<s$, by the Minkowski inequality, applying the H\"{o}lder's
inequality and from (\ref{b}), (\ref{c}), (\ref{314}) we have%
\begin{align*}
\left \Vert J_{4}\right \Vert _{L_{p}\left( E\right) }& \leq \left( \int
\limits_{E}\left \vert \int \limits_{2kr}^{\infty }\int \limits_{E\left(
x_{0},t\right) }\left \vert b\left( y\right) -b_{E\left( x_{0},t\right)
}\right \vert \left \vert f\left( y\right) \right \vert \left \vert \Omega
\left( x-y\right) \right \vert dy\frac{dt}{t^{\gamma +1}}\right \vert
^{p}dx\right) ^{\frac{1}{p}} \\
& +\left( \int \limits_{E}\left \vert \int \limits_{2kr}^{\infty }\left
\vert b_{E\left( x_{0},r\right) }-b_{E\left( x_{0},t\right) }\right \vert
\int \limits_{E\left( x_{0},t\right) }\left \vert f\left( y\right) \right
\vert \left \vert \Omega \left( x-y\right) \right \vert dy\frac{dt}{%
t^{\gamma +1}}\right \vert ^{p}dx\right) ^{\frac{1}{p}} \\
& \lesssim \int \limits_{2kr}^{\infty }\int \limits_{E\left( x_{0},t\right)
}\left \vert b\left( y\right) -b_{E\left( x_{0},t\right) }\right \vert \left
\vert f\left( y\right) \right \vert \left \Vert \Omega \left( \cdot
-y\right) \right \Vert _{L_{p}\left( E\left( x_{0},t\right) \right) }dy\frac{%
dt}{t^{\gamma +1}} \\
& +\int \limits_{2kr}^{\infty }\left \vert b_{E\left( x_{0},r\right)
}-b_{E\left( x_{0},t\right) }\right \vert \int \limits_{E\left(
x_{0},t\right) }\left \vert f\left( y\right) \right \vert \left \Vert \Omega
\left( \cdot -y\right) \right \Vert _{L_{p}\left( E\left( x_{0},t\right)
\right) }dy\frac{dt}{t^{\gamma +1}} \\
& \lesssim \left \vert E\right \vert ^{\frac{1}{p}-\frac{1}{s}}\int
\limits_{2kr}^{\infty }\int \limits_{E\left( x_{0},t\right) }\left \vert
b\left( y\right) -b_{E\left( x_{0},t\right) }\right \vert \left \vert
f\left( y\right) \right \vert \left \Vert \Omega \left( \cdot -y\right)
\right \Vert _{L_{s}\left( E\left( x_{0},t\right) \right) }dy\frac{dt}{%
t^{\gamma +1}} \\
& +\left \vert E\right \vert ^{\frac{1}{p}-\frac{1}{s}}\int
\limits_{2kr}^{\infty }\left \vert b_{E\left( x_{0},r\right) }-b_{E\left(
x_{0},t\right) }\right \vert \int \limits_{E\left( x_{0},t\right) }\left
\vert f\left( y\right) \right \vert \left \Vert \Omega \left( \cdot
-y\right) \right \Vert _{L_{s}\left( E\left( x_{0},t\right) \right) }dy\frac{%
dt}{t^{\gamma +1}} \\
& \lesssim r^{\frac{\gamma }{p}-\frac{\gamma }{s}}\int \limits_{2kr}^{\infty
}\left \Vert \left( b\left( \cdot \right) -b_{E\left( x_{0},t\right)
}\right) \right \Vert _{L_{p_{2}}\left( E\left( x_{0},t\right) \right)
}\left \Vert f\right \Vert _{L_{p_{1}}\left( E\left( x_{0},t\right) \right)
}\left \vert E\left( x_{0},t\right) \right \vert ^{1-\frac{1}{p}}\left \vert
E\left( x_{0},\frac{3}{2}t\right) \right \vert ^{\frac{1}{s}}\frac{dt}{%
t^{\gamma +1}} \\
& +r^{\frac{\gamma }{p}-\frac{\gamma }{s}}\int \limits_{2kr}^{\infty }\left
\vert b_{E\left( x_{0},r\right) }-b_{E\left( x_{0},t\right) }\right \vert
\left \Vert f\right \Vert _{L_{p_{1}}\left( E\left( x_{0},t\right) \right)
}\left \vert E\left( x_{0},\frac{3}{2}t\right) \right \vert ^{\frac{1}{s}}%
\frac{dt}{t^{\frac{\gamma }{p_{1}}+1}} \\
& \lesssim r^{\frac{\gamma }{p}-\frac{\gamma }{s}}\left \Vert b\right \Vert
_{LC_{p_{2},\lambda ,P}^{\left \{ x_{0}\right \} }}\int
\limits_{2kr}^{\infty }\left( 1+\ln \frac{t}{r}\right) t^{\gamma \lambda -%
\frac{\gamma }{p_{1}}+\frac{\gamma }{s}-1}\left \Vert f\right \Vert
_{L_{p_{1}}\left( E\left( x_{0},t\right) \right) }dt.
\end{align*}

Now combined by all the above estimates, we end the proof of this Theorem %
\ref{teo2}.
\end{proof}

Now we can give the following theorem (our main result).

\begin{theorem}
\label{teo9*}Let $x_{0}\in {\mathbb{R}^{n}}$, $1<p<\infty $ and $\Omega \in
L_{s}(S^{n-1})$, $1<s\leq \infty $, be $A_{t}$-homogeneous of degree zero.
Let $T_{\Omega }^{P}$ be a parabolic linear operator satisfying condition (%
\ref{e1}), bounded on $L_{p}({\mathbb{R}^{n}})$ for $1<p<\infty $. Let $b\in
LC_{p_{2},\lambda ,P}^{\left \{ x_{0}\right \} }\left( 
\mathbb{R}
^{n}\right) $, $0\leq \lambda <\frac{1}{\gamma }$ and $\frac{1}{p}=\frac{1}{%
p_{1}}+\frac{1}{p_{2}}$. Let also, for $s^{\prime }\leq p$ the pair $%
(\varphi _{1},\varphi _{2})$ satisfies the condition%
\begin{equation}
\int \limits_{r}^{\infty }\left( 1+\ln \frac{t}{r}\right) \frac{\limfunc{%
essinf}\limits_{t<\tau <\infty }\varphi _{1}(x_{0},\tau )\tau ^{\frac{\gamma 
}{p_{1}}}}{t^{\frac{\gamma }{p_{1}}+1-\gamma \lambda }}dt\leq C\, \varphi
_{2}(x_{0},r),  \label{37}
\end{equation}%
and for $p_{1}<s$ the pair $(\varphi _{1},\varphi _{2})$ satisfies the
condition%
\begin{equation}
\int \limits_{r}^{\infty }\left( 1+\ln \frac{t}{r}\right) \frac{\limfunc{%
essinf}\limits_{t<\tau <\infty }\varphi _{1}(x_{0},\tau )\tau ^{\frac{\gamma 
}{p_{1}}}}{t^{\frac{\gamma }{p_{1}}-\frac{\gamma }{s}+1-\gamma \lambda }}%
dt\leq C\, \varphi _{2}(x_{0},r)r^{\frac{\gamma }{s}},  \label{317*}
\end{equation}%
where $C$ does not depend on $r$.

Then the operator $[b,T_{\Omega }^{P}]$ is bounded from $LM_{p_{1},\varphi
_{1},P}^{\{x_{0}\}}$ to $LM_{p,\varphi _{2},P}^{\{x_{0}\}}$. Moreover,%
\begin{equation}
\left \Vert \lbrack b,T_{\Omega }^{P}]f\right \Vert _{LM_{p,\varphi
_{2},P}^{\{x_{0}\}}}\lesssim \left \Vert b\right \Vert _{LC_{p_{2},\lambda
,P}^{\left \{ x_{0}\right \} }}\left \Vert f\right \Vert _{LM_{p_{1},\varphi
_{1},P}^{\{x_{0}\}}}.  \label{3-1*}
\end{equation}
\end{theorem}

\begin{proof}
Let $p>1$ and $s^{\prime }\leq p$. By Theorem \ref{teo2} and Theorem \ref%
{teo5*} with $v_{2}\left( r\right) =\varphi _{2}\left( x_{0},r\right) ^{-1}$%
, $v_{1}=\varphi _{1}\left( x_{0},r\right) ^{-1}r^{-\frac{\gamma }{p_{1}}}$, 
$w\left( r\right) =r^{\gamma \lambda -\frac{\gamma }{p_{1}}-1}$ and $g\left(
r\right) =\left \Vert f\right \Vert _{L_{p_{1}}\left( E\left( x_{0},r\right)
\right) }$, we have%
\begin{eqnarray*}
\left \Vert \lbrack b,T_{\Omega }^{P}]f\right \Vert _{LM_{p,\varphi
_{2},P}^{\{x_{0}\}}} &\lesssim &\sup_{r>0}\varphi _{2}\left( x_{0},r\right)
^{-1}\left \Vert b\right \Vert _{LC_{p_{2},\lambda ,P}^{\left \{ x_{0}\right \}
}}\int \limits_{r}^{\infty }\left( 1+\ln \frac{t}{r}\right) t^{\gamma \lambda
-\frac{\gamma }{p_{1}}-1}\left \Vert f\right \Vert _{L_{p_{1}}\left( E\left(
x_{0},t\right) \right) }dt \\
&\lesssim &\left \Vert b\right \Vert _{LC_{p_{2},\lambda ,P}^{\left \{
x_{0}\right \} }}\sup_{r>0}\varphi _{1}\left( x_{0},r\right) ^{-1}r^{-\frac{%
\gamma }{p_{1}}}\left \Vert f\right \Vert _{L_{p_{1}}\left( E\left(
x_{0},r\right) \right) }=\left \Vert b\right \Vert _{LC_{p_{2},\lambda
,P}^{\left \{ x_{0}\right \} }}\left \Vert f\right \Vert _{_{LM_{p_{1},\varphi
_{1},P}^{\{x_{0}\}}}},
\end{eqnarray*}%
where the condition (\ref{3}) is equivalent to (\ref{37}), then we obtain (%
\ref{3-1*}).

Let $p>1$ and $p_{1}<s$. By Theorem \ref{teo2} and Theorem \ref{teo5*} with $%
v_{2}\left( r\right) =\varphi _{2}\left( x_{0},r\right) ^{-1}$, $%
v_{1}=\varphi _{1}\left( x_{0},r\right) ^{-1}r^{-\frac{\gamma }{p_{1}}+\frac{%
\gamma }{s}}$, $w\left( r\right) =r^{\gamma \lambda -\frac{\gamma }{p_{1}}+%
\frac{\gamma }{s}-1}$ and $g\left( r\right) =\left \Vert f\right \Vert
_{L_{p_{1}}\left( E\left( x_{0},r\right) \right) }$, we have%
\begin{eqnarray*}
\left \Vert \lbrack b,T_{\Omega }^{P}]f\right \Vert _{_{LM_{p,\varphi
_{2},P}^{\{x_{0}\}}}} &\lesssim &\sup_{r>0}\varphi _{2}\left( x_{0},r\right)
^{-1}r^{-\frac{\gamma }{s}}\left \Vert b\right \Vert _{LC_{p_{2},\lambda
,P}^{\left \{ x_{0}\right \} }}\int \limits_{r}^{\infty }\left( 1+\ln \frac{t}{r%
}\right) t^{\gamma \lambda -\frac{\gamma }{p_{1}}+\frac{\gamma }{s}%
-1}\left \Vert f\right \Vert _{L_{p_{1}}\left( E\left( x_{0},t\right) \right)
}dt \\
&\lesssim &\left \Vert b\right \Vert _{LC_{p_{2},\lambda ,P}^{\left \{
x_{0}\right \} }}\sup_{r>0}\varphi _{1}\left( x_{0},r\right) ^{-1}r^{-\frac{%
\gamma }{p_{1}}}\left \Vert f\right \Vert _{L_{p_{1}}\left( E\left(
x_{0},r\right) \right) }=\left \Vert b\right \Vert _{LC_{p_{2},\lambda
,P}^{\left \{ x_{0}\right \} }}\left \Vert f\right \Vert _{LM_{p_{1},\varphi
_{1},P}^{\{x_{0}\}}},
\end{eqnarray*}%
where the condition (\ref{3}) is equivalent to (\ref{317*}). Thus, we obtain
(\ref{3-1*}).

Hence, the proof is completed.
\end{proof}

In the case of $s=\infty $ from Theorem \ref{teo9*}, we get

\begin{corollary}
Let $x_{0}\in {\mathbb{R}^{n}}$, $1<p<\infty $, $b\in LC_{p_{2},\lambda
,P}^{\left \{ x_{0}\right \} }\left( 
\mathbb{R}
^{n}\right) $, $0\leq \lambda <\frac{1}{\gamma }$, $\frac{1}{p}=\frac{1}{%
p_{1}}+\frac{1}{p_{2}}$ and the pair $(\varphi _{1},\varphi _{2})$ satisfies
condition (\ref{37}). Then the operators $M_{b}^{P}$ and $[b,\overline{T}%
^{P}]$ are bounded from $LM_{p_{1},\varphi _{1},P}^{\{x_{0}\}}$ to $%
LM_{p,\varphi _{2},P}^{\{x_{0}\}}$.
\end{corollary}

\begin{corollary}
Let $x_{0}\in {\mathbb{R}^{n}}$, $1<p<\infty $ and $\Omega \in
L_{s}(S^{n-1}) $, $1<s\leq \infty $, be $A_{t}$-homogeneous of degree zero.
Let $b\in LC_{p_{2},\lambda ,P}^{\left \{ x_{0}\right \} }\left( 
\mathbb{R}
^{n}\right) $, $0\leq \lambda <\frac{1}{\gamma }$, $\frac{1}{p}=\frac{1}{%
p_{1}}+\frac{1}{p_{2}}$. Let also, for $s^{\prime }\leq p$, the pair $\left(
\varphi _{1},\varphi _{2}\right) $ satisfies condition (\ref{37}) and for $%
p<s$, the pair $\left( \varphi _{1},\varphi _{2}\right) $ satisfies
condition (\ref{317*}). Then the operators $M_{\Omega ,b}^{P}$ and $[b,%
\overline{T}_{\Omega }^{P}]$ are bounded from $LM_{p_{1},\varphi
_{1},P}^{\{x_{0}\}}$ to $LM_{p,\varphi _{2},P}^{\{x_{0}\}}$.
\end{corollary}

\begin{corollary}
\label{Corollary 1*}Let $x_{0}\in {\mathbb{R}^{n}}$, $1<p<\infty $ and $%
\Omega \in L_{s}(S^{n-1})$, $1<s\leq \infty $, be homogeneous of degree
zero. Let $T_{\Omega }^{P}$ be a parabolic linear operator satisfying
condition (\ref{e1}), bounded on $L_{p}({\mathbb{R}^{n}})$ for $1<p<\infty $%
. Let $b\in LC_{p_{2},\lambda ,P}^{\left \{ x_{0}\right \} }\left( 
\mathbb{R}
^{n}\right) $, $0\leq \lambda <\frac{1}{n}$, $\frac{1}{p}=\frac{1}{p_{1}}+%
\frac{1}{p_{2}}$. Let also, for $s^{\prime }\leq p$ the pair $(\varphi
_{1},\varphi _{2})$ satisfies the condition%
\begin{equation}
\int \limits_{r}^{\infty }\left( 1+\ln \frac{t}{r}\right) \frac{\limfunc{%
essinf}\limits_{t<\tau <\infty }\varphi _{1}(x_{0},\tau )\tau ^{\frac{n}{%
p_{1}}}}{t^{\frac{n}{p_{1}}+1-n\lambda }}dt\leq C\, \varphi _{2}(x_{0},r),
\label{ferit}
\end{equation}%
and for $p_{1}<s$ the pair $(\varphi _{1},\varphi _{2})$ satisfies the
condition%
\begin{equation*}
\int \limits_{r}^{\infty }\left( 1+\ln \frac{t}{r}\right) \frac{\limfunc{%
essinf}\limits_{t<\tau <\infty }\varphi _{1}(x_{0},\tau )\tau ^{\frac{n}{%
p_{1}}}}{t^{\frac{n}{p_{1}}-\frac{n}{s}+1-n\lambda }}dt\leq C\, \varphi
_{2}(x_{0},r)r^{\frac{n}{s}},
\end{equation*}%
where $C$ does not depend on $r$.

Then the operator $[b,T_{\Omega }^{P}]$ is bounded from $LM_{p_{1},\varphi
_{1},P}^{\{x_{0}\}}$ to $LM_{p,\varphi _{2},P}^{\{x_{0}\}}$. Moreover,%
\begin{equation*}
\left \Vert \lbrack b,T_{\Omega }^{P}]f\right \Vert _{LM_{p,\varphi
_{2},P}^{\{x_{0}\}}}\lesssim \left \Vert b\right \Vert _{LC_{p_{2},\lambda
,P}^{\left \{ x_{0}\right \} }}\left \Vert f\right \Vert _{LM_{p_{1},\varphi
_{1},P}^{\{x_{0}\}}}.
\end{equation*}
\end{corollary}

\begin{remark}
Note that, in the case of $P=I$ Corollary \ref{Corollary 1*} has been proved
in \cite{BGGS, Gurbuz1}. Also, in the case of $P=I$ and $s=\infty $
Corollary \ref{Corollary 1*} has been proved in \cite{BGGS, Gurbuz1}.
\end{remark}

\begin{corollary}
\label{Corollary 2*}Let $\Omega \in L_{s}(S^{n-1})$, $1<s\leq \infty $, be $%
A_{t}$-homogeneous of degree zero. Let $T_{\Omega }^{P}$ be a parabolic
linear operator satisfying condition (\ref{e1}), bounded on $L_{p}({\mathbb{R%
}^{n}})$ for $1<p<\infty $. Let $1<p<\infty $ and $b\in BMO\left( 
\mathbb{R}
^{n}\right) $ (parabolic bounded mean oscillation space). Let also, for $%
s^{\prime }\leq p$ the pair $(\varphi _{1},\varphi _{2})$ satisfies the
condition%
\begin{equation}
\int \limits_{r}^{\infty }\left( 1+\ln \frac{t}{r}\right) \frac{\limfunc{%
essinf}\limits_{t<\tau <\infty }\varphi _{1}(x,\tau )\tau ^{\frac{\gamma }{p}%
}}{t^{\frac{\gamma }{p}+1}}dt\leq C\, \varphi _{2}(x,r),  \label{37*}
\end{equation}%
and for $p<s$ the pair $(\varphi _{1},\varphi _{2})$ satisfies the condition%
\begin{equation}
\int \limits_{r}^{\infty }\left( 1+\ln \frac{t}{r}\right) \frac{\limfunc{%
essinf}\limits_{t<\tau <\infty }\varphi _{1}(x,\tau )\tau ^{\frac{\gamma }{p}%
}}{t^{\frac{\gamma }{p}-\frac{\gamma }{s}+1}}dt\leq C\, \varphi _{2}(x,r)r^{%
\frac{\gamma }{s}},  \label{1}
\end{equation}%
where $C$ does not depend on $x$ and $r$.

Then the operator $[b,T_{\Omega }^{P}]$ is bounded from $M_{p,\varphi
_{1},P} $ to $M_{p,\varphi _{2},P}$. Moreover,%
\begin{equation*}
\left \Vert \lbrack b,T_{\Omega }^{P}]f\right \Vert _{M_{p,\varphi
_{2},P}}\lesssim \left \Vert b\right \Vert _{BMO}\left \Vert f\right \Vert
_{M_{p,\varphi _{1},P}}.
\end{equation*}
\end{corollary}

In the case of $s=\infty $ from Corollary \ref{Corollary 2*}, we get

\begin{corollary}
Let $1<p<\infty $, $b\in BMO\left( 
\mathbb{R}
^{n}\right) $ and the pair $(\varphi _{1},\varphi _{2})$ satisfies condition
(\ref{37*}). Then the operators $M_{b}^{P}$ and $[b,\overline{T}^{P}]$ are
bounded from $M_{p,\varphi _{1},P}$ to $M_{p,\varphi _{2},P}$.
\end{corollary}

\begin{corollary}
Let $\Omega \in L_{s}(S^{n-1})$, $1<s\leq \infty $, be $A_{t}$-homogeneous
of degree zero. Let $1<p<\infty $ and $b\in BMO\left( 
\mathbb{R}
^{n}\right) $. Let also, for $s^{\prime }\leq p$, the pair $\left( \varphi
_{1},\varphi _{2}\right) $ satisfies condition (\ref{37*}) and for $p<s$,
the pair $\left( \varphi _{1},\varphi _{2}\right) $ satisfies condition (\ref%
{1}). Then the operators $M_{\Omega ,b}^{P}$ and $[b,\overline{T}_{\Omega
}^{P}]$ are bounded from $M_{p,\varphi _{1},P}$ to $M_{p,\varphi _{2},P}$.
\end{corollary}

\begin{corollary}
\label{Corollary 3*}Let $\Omega \in L_{s}(S^{n-1})$, $1<s\leq \infty $, be
homogeneous of degree zero. Let $1<p<\infty $ and $b\in BMO\left( 
\mathbb{R}
^{n}\right) $. Let also, for $s^{\prime }\leq p$ the pair $(\varphi
_{1},\varphi _{2})$ satisfies the condition%
\begin{equation*}
\int \limits_{r}^{\infty }\left( 1+\ln \frac{t}{r}\right) \frac{\limfunc{%
essinf}\limits_{t<\tau <\infty }\varphi _{1}(x,\tau )\tau ^{\frac{n}{p}}}{t^{%
\frac{n}{p}+1}}dt\leq C\, \varphi _{2}(x,r),
\end{equation*}%
and for $p<s$ the pair $(\varphi _{1},\varphi _{2})$ satisfies the condition%
\begin{equation*}
\int \limits_{r}^{\infty }\left( 1+\ln \frac{t}{r}\right) \frac{\limfunc{%
essinf}\limits_{t<\tau <\infty }\varphi _{1}(x,\tau )\tau ^{\frac{n}{p}}}{t^{%
\frac{n}{p}-\frac{n}{s}+1}}dt\leq C\, \varphi _{2}(x,r)r^{\frac{n}{s}},
\end{equation*}%
where $C$ does not depend on $x$ and $r$.

Then the operator $[b,T_{\Omega }^{P}]$ is bounded from $M_{p,\varphi
_{1},P} $ to $M_{p,\varphi _{2},P}$. Moreover,%
\begin{equation*}
\left \Vert \lbrack b,T_{\Omega }^{P}]f\right \Vert _{M_{p,\varphi
_{2},P}}\lesssim \left \Vert b\right \Vert _{BMO}\left \Vert f\right \Vert
_{M_{p,\varphi _{1},P}}.
\end{equation*}
\end{corollary}

\begin{remark}
Note that, in the case of $P=I$ Corollary \ref{Corollary 3*} has been proved
in \cite{BGGS, Gurbuz1, Gurbuz2}. Also, in the case of $P=I$ and $s=\infty $
Corollary \ref{Corollary 3*} has been proved in \cite{BGGS, GULAKShIEOT2012,
Gurbuz1, Gurbuz2, Karaman}.
\end{remark}

Now, we give the applications of Theorem \ref{teo9} and Theorem \ref{teo9*}
for the parabolic Marcinkiewicz operator.

Suppose that $\Omega \left( x\right) $ is a real-valued and measurable
function defined on ${\mathbb{R}^{n}}$ satisfying the following conditions:

(a) $\Omega \left( x\right) $ is homogeneous of degree zero with respect to $%
A_{t}$, that is, 
\begin{equation*}
\Omega (A_{t}x)=\Omega (x),~\text{for any}~~t>0,\text{ }x\in {\mathbb{R}^{n}}%
\setminus \{0\};
\end{equation*}

(b) $\Omega \left( x\right) $ has mean zero on $S^{n-1}$, that is, 
\begin{equation*}
\dint \limits_{S^{n-1}}\Omega (x^{\prime })J\left( x^{\prime }\right)
d\sigma (x^{\prime })=0,
\end{equation*}%
where $x^{\prime }=\frac{x}{\left \vert x\right \vert }$ for any $x\neq 0$.

(c) $\Omega \in L_{1}\left( S^{n-1}\right) $.

Then the parabolic Marcinkiewicz integral of higher dimension $\mu _{\Omega
}^{\gamma }$ is defined by 
\begin{equation*}
\mu _{\Omega }^{\gamma }(f)(x)=\left( \int \limits_{0}^{\infty }|F_{\Omega
,t}(f)(x)|^{2}\frac{dt}{t^{3}}\right) ^{1/2},
\end{equation*}%
where 
\begin{equation*}
F_{\Omega ,t}(f)(x)=\int \limits_{{\rho }\left( x-y\right) \leq t}\frac{%
\Omega (x-y)}{{\rho }\left( x-y\right) ^{\gamma -1}}f(y)dy.
\end{equation*}

On the other hand, for a suitable function $b$, the commutator of the
parabolic Marcinkiewicz integral $\mu _{\Omega }^{\gamma }$ is defined by 
\begin{equation*}
\lbrack b,\mu _{\Omega }^{\gamma }](f)(x)=\left( \dint \limits_{0}^{\infty
}|F_{\Omega ,t,b}(f)(x)|^{2}\frac{dt}{t^{3}}\right) ^{1/2},
\end{equation*}%
where 
\begin{equation*}
F_{\Omega ,t,b}(f)(x)=\dint \limits_{{\rho }\left( x-y\right) \leq t}\frac{%
\Omega (x-y)}{{\rho }\left( x-y\right) ^{\gamma -1}}[b(x)-b(y)]f(y)dy.
\end{equation*}

We consider the space $H=\{h:\Vert h\Vert =(\int \limits_{0}^{\infty
}|h(t)|^{2}\frac{dt}{t^{3}})^{1/2}<\infty \}$. Then, it is clear that $\mu
_{\Omega }^{\gamma }(f)(x)=\Vert F_{\Omega ,t}(x)\Vert $.

By the Minkowski inequality and the conditions on $\Omega $, we get 
\begin{equation*}
\mu _{\Omega }^{\gamma }(f)(x)\leq \int \limits_{{\mathbb{R}^{n}}}\frac{%
|\Omega (x-y)|}{{\rho }\left( x-y\right) ^{\gamma -1}}|f(y)|\left( \int
\limits_{|x-y|}^{\infty }\frac{dt}{t^{3}}\right) ^{1/2}dy\leq C\int \limits_{%
{\mathbb{R}^{n}}}\frac{\left \vert \Omega (x-y)\right \vert }{{\rho }\left(
x-y\right) ^{\gamma }}|f(y)|dy.
\end{equation*}%
Thus, $\mu _{\Omega }^{\gamma }$ satisfies the condition (\ref{e1}). When $%
\Omega \in L_{s}\left( S^{n-1}\right) \left( s>1\right) $, It is known that $%
\mu _{\Omega }$ is bounded on $L_{p}({\mathbb{R}^{n}})$ for $p>1$, and
bounded from $L_{1}({\mathbb{R}^{n}})$ to $WL_{1}({\mathbb{R}^{n}})$ for $%
p=1 $ (see \cite{XDY}), then from Theorems \ref{teo9}, \ref{teo9*} we get

\begin{corollary}
Let $x_{0}\in {\mathbb{R}^{n}}$, $1\leq p<\infty $, $\Omega \in L_{s}\left(
S^{n-1}\right) $, $1<s\leq \infty $. Let also, for $s^{\prime }\leq p$, $%
p\neq 1$, the pair $\left( \varphi _{1},\varphi _{2}\right) $ satisfies
condition (\ref{e37}) and for $1<$ $p<s$ the pair $\left( \varphi
_{1},\varphi _{2}\right) $ satisfies condition (\ref{317}) and $\Omega $
satisfies conditions (a)--(c). Then the operator $\mu _{\Omega }^{\gamma }$
is bounded from $LM_{p,\varphi _{1}}^{\{x_{0}\}}$ to $LM_{p,\varphi
_{2}}^{\{x_{0}\}}$ for $p>1$ and from $LM_{1,\varphi _{1}}^{\{x_{0}\}}$ to $%
WLM_{1,\varphi _{2}}^{\{x_{0}\}}$.
\end{corollary}

\begin{corollary}
Let $1\leq p<\infty $, $\Omega \in L_{s}\left( S^{n-1}\right) $, $1<s\leq
\infty $. Let also, for $s^{\prime }\leq p$, $p\neq 1$, the pair $\left(
\varphi _{1},\varphi _{2}\right) $ satisfies condition (\ref{e37*}) and for $%
1<$ $p<s$ the pair $\left( \varphi _{1},\varphi _{2}\right) $ satisfies
condition (\ref{e38}) and $\Omega $ satisfies conditions (a)--(c). Then the
operator $\mu _{\Omega }^{\gamma }$ is bounded from $M_{p,\varphi _{1}}$ to $%
M_{p,\varphi _{2}}$ for $p>1$ and from $M_{1,\varphi _{1}}$to $WM_{1,\varphi
_{2}}$ for $p=1$.
\end{corollary}

\begin{corollary}
Let $x_{0}\in {\mathbb{R}^{n}}$, $\Omega \in L_{s}(S^{n-1})$, $1<s\leq
\infty $. Let $1<p<\infty $, $b\in LC_{p_{2},\lambda }^{\left \{
x_{0}\right
\} }({\mathbb{R}^{n}})$, $\frac{1}{p}=\frac{1}{p_{1}}+\frac{1}{%
p_{2}}$, $0\leq \lambda <\frac{1}{n}$. Let also, for $s^{\prime }\leq p$ the
pair $(\varphi _{1},\varphi _{2})$ satisfies condition (\ref{37}) and for $%
p_{1}<s$ the pair $(\varphi _{1},\varphi _{2})$ satisfies condition (\ref%
{317*}) and $\Omega $ satisfies conditions (a)--(c). Then, the operator $%
[b,\mu _{\Omega }^{\gamma }]$ is bounded from $LM_{p_{1},\varphi
_{1}}^{\{x_{0}\}}$ to $LM_{p,\varphi _{2}}^{\{x_{0}\}}$.
\end{corollary}

\begin{corollary}
Let $\Omega \in L_{s}(S^{n-1})$, $1<s\leq \infty $, $1<p<\infty $ and $b\in
BMO\left( 
\mathbb{R}
^{n}\right) $. Let also, for $s^{\prime }\leq p$ the pair $(\varphi
_{1},\varphi _{2})$ satisfies condition (\ref{37*}) and for $p<s$ the pair $%
(\varphi _{1},\varphi _{2})$ satisfies condition (\ref{1}) and $\Omega $
satisfies conditions (a)--(c). Then, the operator $[b,\mu _{\Omega }^{\gamma
}]$ is bounded from $M_{p,\varphi _{1}}$ to $M_{p,\varphi _{2}}$.
\end{corollary}

\begin{remark}
Obviously, If we take $\alpha _{1}=\cdots \alpha _{n}=1$ and $P=I$, then
obviously $\rho \left( x\right) =\left \vert x\right \vert =\left( \dsum
\limits_{i=1}^{n}x_{i}^{2}\right) ^{\frac{1}{2}}$, $\gamma =n$, $\left( {%
\mathbb{R}^{n},\rho }\right) =$ $\left( {\mathbb{R}^{n},}\left \vert \cdot
\right \vert \right) $, $E_{I}(x,r)=B\left( x,r\right) $. In this case, $\mu
_{\Omega }^{\gamma }$ is just the classical Marcinkiewicz integral operator $%
\mu _{\Omega }$, which was first defined by Stein in 1958. In \cite{Stein58}%
, Stein has proved that if $\Omega $ satisfies the Lipshitz condition of
degree of $\alpha \left( 0<\alpha \leq 1\right) $ on $S^{n-1}$ and the
conditions (a), (b) (obviously, in the case $A_{t}=tI$ and $J\left(
x^{\prime }\right) \equiv 1$), then $\mu _{\Omega }$ is both of the type $%
\left( p,p\right) $ $\left( 1<p\leq 2\right) $ and the weak type $\left(
1.1\right) $. (See also \cite{TorWang} for the boundedness of the classical
Marcinkiewicz integral $\mu _{\Omega }$.)
\end{remark}

\end{document}